\documentclass[12pt,reqno]{amsart}

\usepackage{amssymb,latexsym}
\usepackage{graphicx,psfrag,color}
\usepackage[square]{natbib}

\newtheorem{thm}{Theorem}[section]
\newtheorem{prop}[thm]{Proposition}

\newtheorem{lem}[thm]{Lemma}
\newtheorem{dfn}[thm]{Definition}
\newtheorem{exa}[thm]{Example}
\theoremstyle{definition}
\newtheorem{rem}[thm]{Remark}

\allowdisplaybreaks[1]

\makeatletter\@addtoreset{equation}{section}\makeatother

\def\bk{\rr{K\B}}
\def\ct{\mathop{\rm CT}}
\def\CT{\mathop{\rm CT}}
\def\pa{\mathbb{ A}_s}
\def\pb{\mathbb{B}_s}
\def\trpath{$\rr{\B}$-path }
\def\krpath{$K_r$-path }
\def\ktrpath{$\rr{K\B}$-path }

\def\C{\hat{C}}
\def\B{T}
\makeatletter
\def\asterisk{\ifnum\c@page=1 $^*$ \fi}
\makeatother

\newcommand{\rr}[1]{#1^{(r)}}

\newcommand{\thyper}[3]{{}_2F_1 \left(#1; #2\; \bigg|\; #3\right)}
\newcommand{\shyper}[3]{{}_2F_1 \left(#1; #2\mid #3\right)}

\newcommand{\tf}{\frac{27}4}

\renewcommand{\r}[2]{(#1)_{#2}}

\def\ijn{_{0\le i,j\le n-1}}
\newcommand{\detijn}[1]{\det\left(#1\right)\ijn}
\begin{document}
\thanks{$^*$Partially supported by NSF Grant DMS-0200596.}

\title[The Generating Function of Ternary
Trees and Continued Fractions]{The Generating Function of Ternary
Trees and Continued Fractions}
\author{Ira M. Gessel\asterisk}
\address{Department of
Mathematics\\
Brandeis University\\
Waltham MA
02454-9110}
\email{gessel@brandeis.edu}
\author{Guoce Xin}
\address{Department
of Mathematics\\
Brandeis University\\
Waltham MA
02454-9110}
\email{guoce.xin@gmail.com}
\date{May 11, 2005}
\begin{abstract}
Michael Somos conjectured a relation between Hankel
determinants whose entries  $\frac 1{2n+1}\binom{3n}n$
count ternary trees and the number of
certain plane partitions and  alternating sign
matrices. Tamm evaluated these determinants by showing that the
generating function
for these entries has a continued fraction that is a special case
of Gauss's continued fraction for a quotient of hypergeometric
series. We give a systematic application of the continued fraction
method to a number of similar Hankel determinants. We also describe a
simple method for transforming determinants using the generating
function for their entries. In this way we transform
Somos's Hankel determinants to known determinants,  and we obtain, up to a power of $3$, a Hankel determinant for the number
of alternating sign matrices. We obtain a combinatorial proof, in
terms of nonintersecting paths, of determinant identities
involving the number of ternary trees and more general determinant identities involving the number of $r$-ary trees.
\end{abstract}
\maketitle


\vspace{3mm} \noindent {\small {\em \small Mathematics Subject
Classifications}: Primary 05A15, Secondary 05A10, 05A17, 30B70,
33C05}

\section{Introduction}

Let $a_n=\frac1{2n+1}\binom{3n}{n}=\frac1{3n+1}\binom{3n+1}n$ be the number of
ternary trees with $n$ vertices and define the Hankel determinants
\begin{align}
U_n&=\detijn {a_{i+j}} \label{e-s1} \\
V_n&=\detijn{a_{i+j+1}} \label{e-s2}\\
W_n&=\detijn{a_{(i+j+1)/2}}, \label{e-s3}
\end{align}
where we take $a_k$ to be 0 if $k$ is not an integer. (We also interpret determinants of $0\times 0$
matrices as 1.) The first few values of these determinants are
$$
\begin{array}{r|rrrrrrr}
n&1\hfil&2\hfil&3\hfil&4\hfil&5\hfil&6\hfil&7\hfil\\  \hline
U_n&1&2&11&170&7429&920460&323801820\\
V_n&1&3&26&646&45885&9304650&5382618660\\
W_n&1&1&2&6&33&286&4420
\end{array}$$

This paper began as an attempt to prove the conjectures of Michael Somos \cite{somos} that
\begin{enumerate}
\item[(a)] $U_n$ is the number of of cyclically symmetric transpose complement plane partitions whose Ferrers diagrams
fit in an $n\times n\times n$ box,
\item[(b)] $V_n$ is the number of $(2n+1)\times (2n+1)$ alternating sign matrices that are invariant under vertical
reflection, and
\item[(c)] $W_n$ is the number of $(2n+1)\times (2n+1)$
alternating sign  matrices that are invariant under both vertical and horizontal reflection.
\end{enumerate}

Mills, Robbins, and Rumsey \cite{tcsymm} (see also \cite[Eq.
(6.15), p. 199]{paproof}) showed  that the number of objects of
type (a) is
\begin{align}
\prod_{i=1}^{n-1} \frac{(3i+1)(6i)!\,(2i)!}{(4i+1)!\,(4i)!}.
\label{e-typea}
\end{align}
Mills  \cite{robbins}
conjectured the formula
\begin{equation}
\prod_{i=1}^{n}\frac{\binom {6i-2}{2i}}{2\binom{4i-1}{2i}}
\label{e-typeb}
\end{equation}
for objects of type (b)
and this conjecture was proved by Kuperberg \cite{kuperberg}. A formula for objects of type (c)
was conjectured by Robbins \cite{robbins2} and proved by
Okada \cite{okada}. A determinant
formula for these objects was proved by Kuperberg \cite{kuperberg}.

It turns out that it is much easier to evaluate Somos's determinants than to relate them to (a)--(c).  It is easy to see that
$W_{2n}=U_nV_n$ and
$W_{2n+1}=U_{n+1}V_n$, so it is only
necessary show that  $U_n$ is equal to \eqref{e-typea} and
$V_n$ is equal to \eqref{e-typeb} to prove Somos's conjectures.

This was done by Tamm \cite{tamm}, who was unaware of Somos's
conjectures. Thus Somos's conjectures are already proved; nevertheless, our study of these conjecture led to some additional determinant evaluations and transformations that are the subject of this paper.

Tamm's proof used the fact that Hankel determinants can be
evaluated using continued fractions; the continued fraction that
gives these Hankel determinants is a special case of Gauss's
continued fraction for a quotient of hypergeometric series.
The determinant $V_n$ was also evaluated, using a different
method, by E\u gecio\u glu, Redmond, and Ryavec \cite[Theorem 4]{err}, who also noted the connection with alternating sign matrices and gave several additional
Hankel determinants for $V_n$:
\begin{equation}
\label{e-err}
V_n=\detijn{b_{i+j}}=\detijn{r_{i+j}}=\detijn {s_{i+j}(u)},
\end{equation}
where $b_n=\frac1{n+1}\binom {3n+1}n$, $r_n=\binom{3n+2}n$, and
\[s_n(u)=\sum_{k=0}^{n}\frac{k+1}{n+1}\binom{3n-k+1}{n-k}u^k,\]
where $u$ is arbitrary. As noted in \cite[Theorem~4]{err},
$s_n(0)=b_n$, $s_n(1)=a_{n+1}$, and $s_n(3)=r_n$.

In Section 2, we describe Tamm's continued fraction method for evaluating these determinants.
 In Section 3, we give a systematic application of the continued fraction method to several similar Hankel determinants.
In Theorem \ref{t-cfs} we
give  five pairs  of generating functions similar to that for $a_n$
whose continued fractions are instances of Gauss's theorem.
Three of them have known
combinatorial meanings for their coefficients, including the number
of two-stack-sortable permutations (see West \cite{west}).

In Section 4 we discuss a simple method, using generating functions, for transforming determinants and use it to show that
\begin{align}
U_n&=\detijn{\binom{i+j}{2i-j}}\label{e-un} \\
\intertext{and} V_n&=\detijn{\binom{i+j+1}{2i-j}}. \label{e-vn}
\end{align}

We also prove E\u gecio\u glu, Redmond, and Ryavec's identity
\eqref{e-err} and the related identity
\begin{align}
\detijn{s_{i+j-1}(u)}&= U_n/u,\quad n>0, \label{e-sdet}
\end{align}
where $s_{-1}(u)=u^{-1}$. When $u=1$, \eqref{e-sdet} reduces to
\eqref{e-s1} and  when $u=3$, \eqref{e-sdet} reduces to
\begin{equation}
\label{e-rdet}
\detijn{r_{i+j-1}}=U_{n}/3, \quad n>0.
\end{equation}
Note that $r_{n-1}=\frac13\binom{3n}n$, so \eqref{e-rdet} is
equivalent to $\detijn{\binom{3n}n}=3^{n-1}U_n$ for $n>0$.

In Section 5 we consider the Hankel determinants of the
coefficients of
\[\frac{1-(1-9x)^{1/3}}{3x}.\]
We first evaluate them using continued fractions, and then show
that the method of Section 4 transforms them into powers of 3
times the determinant
\[\detijn{\binom{i+j}{i-1}+\delta_{ij}},\]
which counts descending plane partitions and alternating sign matrices.
Similarly, the Hankel determinant corresponding
to
\[\frac{1-(1-9x)^{2/3}}{3x}\]
is transformed to a power of 3 times the determinant
\[\detijn{\binom{i+j}{i}+\delta_{ij}},\]
which counts cyclically symmetric plane partitions.

By a result of Gessel and Viennot \cite{gessel-viennot}, both
sides of \eqref{e-un} can be interpreted by $n$-tuples of
nonintersecting lattice paths. A similar situation holds for
\eqref{e-vn}. In Section 6, we describe the nonintersecting
lattice path interpretations of these determinants. We give a new
class of interpretations of $a_n$ in terms of certain paths called
$K$-paths in Theorem \ref{t-newg}. {}From this new interpretation
of $a_n$, \eqref{e-un} follows easily. The proof of Theorem
\ref{t-newg} relies on a ``sliding lemma," which says that the
number of certain $K$-paths does not change after sliding their
starting and ending points.

In Section 7, we study another class of paths called $T$-paths, which are related to trinomial
coefficients, and $KT$-paths, which are analogous to $K$-paths. We find another class of interpretations of $a_n$ in
terms of $KT$-paths, using which we find a new determinant
identity involving $U_n$ (Theorem \ref{t-kte}). Unfortunately,
we do not have a nonintersecting path interpretation for this determinant. There is a
natural bijection from $K$-paths to $KT$-paths, and the sliding
lemma for $KT$-paths is easier to prove than that for $K$-paths.

In Section 8, we study  $KT^{(r)}$-paths, which reduce to $KT$ paths when $r=2$. The results of Section~7 generalize, and  we obtain determinant identities involving Hankel
determinants for the number of
$(r+1)$-ary trees  (see \eqref{e-98} and \eqref{e-99}).

In Section 9, we give algebraic proofs of the results of Section 8
using  partial fractions.

\section{Hankel Determinants and Gauss's Continued Fraction \label{sec-Hankel-Gauss}}
Let $A(x)=\sum_{n\geq 0} A_n x^n$ be a formal power series. We
define the Hankel determinants $H_n^{(k)}(A)$  of $A(x)$ by
\begin{equation*}
H_n^{(k)}(A)=\detijn{A_{i+j+k}}.
\end{equation*}
We shall write $H_n(A)$ for $H_n^{(0)}(A)$ and $H_n^1(A)$ for
$H_n^{(1)}(A)$. We also define ${\hat H}_n(A)$ to be
$H_n(A(x^2))$.
 It is not difficult to show that
${\hat H}_{2n}(A)=H_n(A)H^1_n(A)$ and
${\hat H}_{2n+1}(A)=H_{n+1}(A)H^1_n(A)$.

Let $g(x)$ be the generating function for ternary trees:
\begin{equation}\label{e-g}
g(x)=\sum_{n\geq 0} a_n  x^n=\sum_{n\geq 0}
\frac{1}{2n+1}{3n\choose n} x^n,
\end{equation}
which is uniquely determined by the functional equation
\begin{equation}\label{e-gf}
g(x)=1+x g(x)^3.
\end{equation}
Then $U_n=H_n(g(x))$, $V_n=H_n^1(g(x))$, and
$W_n=\hat{H}_n(g(x))$.

In general, it is difficult to say much about $H_n(A(x))$.
However, if $A(x)$ can be expressed as a continued fraction, then there
is a very nice formula. This is the case for $g(x)$:
Tamm \cite{tamm} observed that
$g(x)$ has a nice continued fraction
expression, which is a special case of Gauss's continued fraction.
We introduce some notation to explain Tamm's approach.

We use the notation $S(x;\lambda_1,\lambda_2,\lambda_3,\ldots )$
to denote the continued fraction
\begin{equation}
S(x;\lambda_1,\lambda_2,\lambda_3,\ldots
)=
\displaystyle\frac{1}{1-\displaystyle\frac{
\lambda_1
x}{1-\displaystyle\frac{\lambda_2 x}{1-\displaystyle\frac{\lambda_3
x}{\ddots}}} }
\end{equation}
The following theorem is equivalent to
\cite[Theorem 7.2]{cfraction}.
Additional information about continued fractions and Hankel determinants can be found in Krattenthaler \cite[Section 5.4]{kratt2}.

\begin{lem}
\label{l-Hcf}
Let
$A(x)=S(x;\lambda_1,\lambda_2,\lambda_3,\ldots )$ and let
$\mu_i=\lambda_1\lambda_2\cdots \lambda_i$.
Then for $n\ge1$,
\begin{align}
H_n(A)
& =  (\lambda_1\lambda_2)^{n-1} (\lambda_3\lambda_4)^{n-2}
    \cdots (\lambda_{2n-3}\lambda_{2n-2})
    =\mu_2\mu_4\cdots \mu_{2n-2}  \label{e-hna}\\
H_n^1(A) & =  \lambda_1^n (\lambda_2\lambda_3)^{n-1}
  \cdots (\lambda_{2n-2}\lambda_{2n-1})
   =\mu_1\mu_3\cdots\mu_{2n-1} \label{e-kna}\\
{\hat H}_n(A) & = \lambda_1^{n-1}\lambda_2^{n-2}\cdots
\lambda_{n-2}^2\lambda_{n-1}=\mu_1\mu_2\cdots \mu_{n-1}.\label{3-hata}
\end{align}
\end{lem}

We define the hypergeometric series by
\[\shyper{a,b}cx=\sum_{n=0}^\infty \frac{\r an \r bn}{n!\, \r
cn} x^n,\] where $\r un=u(u+1)\cdots (u+n-1)$.

Gauss proved the following
theorem \cite[Theorem 6.1]{cfraction}, which gives a continued
fraction for a quotient of two hypergeometric series:
\begin{lem}\label{l-gauss}
If $c$ is not a negative integer then we have the continued fraction
\begin{equation}
\shyper{a,b+1}{c+1}{x}
\big /\shyper{a,b}{c}{x}
=S(x;\lambda_1,\lambda_2,\ldots ),
\label{e-gausscf}
\end{equation}
where
\begin{equation}\label{e-lamg}
\begin{aligned}
\lambda_{2n-1} & =\dfrac{(a+n-1)(c-b+n-1)}{(c+2n-2)(c+2n-1)},
&\quad& n=1,2,\ldots,
\\
\lambda_{2n} & = \dfrac{(b+n)(c-a+n)}{(c+2n-1)(c+2n)},
&&n=1,2,\ldots .
\end{aligned}
\end{equation}
\end{lem}

Combining Lemmas \ref{l-Hcf}  and \ref{l-gauss} gives a formula for evaluating certain Hankel
determinants.

\begin{lem}
\label{l-hankel}
Let
\[A(x)=\shyper{a,b+1}{c+1}{\rho x}\big/ \shyper{a,b}c{\rho x}.\]
Then
\begin{align}
H_n(A)&=\prod_{i=0}^{n-1}\frac {\r ai \r {b+1}i \r{c-b}i  \r{c-a+1}i}
  {\r c{2i} \r {c+1}{2i}}  \rho^{2i} \label{e-h1}\\
H_n^1(A)&=\prod_{i=1}^{n}\frac {\r ai \r {b+1}{i-1} \r{c-b}i  \r{c-a+1}{i-1}}
  {\r c{2i-1} \r {c+1}{2i-1}}  \rho^{2i-1} \label{e-h2}\\
        &=\prod_{i=1}^{n}\frac{(c-1)c} {b(c-a)\rho}
       \frac {\r ai \r {b}{i} \r{c-b}i  \r{c-a}{i}}
        {\r c{2i} \r {c-1}{2i}} \rho^{2i} \label{e-h3}
\end{align}
\end{lem}

\begin{proof}
By Lemma \ref{l-gauss}, $A(x)$ has the continued fraction expansion
$A(x)=S(x; \lambda_1, \lambda_2,\cdots)$
where
\begin{align*}
\lambda_{2n-1} & =\frac{(a+n-1)(c-b+n-1)}{(c+2n-2)(c+2n-1)}\rho,
\\
\lambda_{2n} & = \frac{(b+n)(c-a+n)}{(c+2n-1)(c+2n)}\rho.
\end{align*}
Then
\begin{align*}
&\lambda_1\lambda_3\cdots\lambda_{2i-1}=\frac{\r ai\r{c-b} i}{\r c{2i}}\rho^i \\
\intertext{and}
&\lambda_2\lambda_4\cdots \lambda_{2i}=\frac{\r {b+1}i \r{c-a+1}i}{\r{c+1}{2i}}\rho^i.
\end{align*}
So with the notation of Lemma \ref{l-Hcf},
\begin{align*}
\mu_{2i}&=\lambda_1\lambda_2\cdots \lambda_{2i}=\frac {\r ai \r{c-b}i \r {b+1}i \r{c-a+1}i}
  {\r c{2i} \r {c+1}{2i}}  \rho^{2i}  \\
\intertext{and}
\mu_{2i-1}&=\lambda_1\lambda_2\cdots \lambda_{2i-1}
=\frac {\r ai \r{c-b}i \r {b+1}{i-1} \r{c-a+1}{i-1}}
  {\r c{2i} \r {c+1}{2i-2}}  \rho^{2i-1}.
\end{align*}
Then \eqref{e-h1} follows immediately from \eqref{e-hna}, and \eqref{e-h2} follows from
\eqref{e-kna} with the help of the identity $\r c{2i} \r {c+1}{2i-2}= \r c{2i-1} \r{c+1}{2i-1}$, and
\eqref{e-h3} follows easily from \ref{e-h2}.
\end{proof}

There is also a simple formula for $H^{(2)}_n(A)$, although we will not need it.

\begin{lem}
\label{l-Q}
 Let $Q(a,b,c\mid x)=
\shyper{a,b+1}{c+1}x / \shyper{a,b}cx$.
Then
\[Q(b,a,c\mid x)=\frac{c(a-b)}{a(c-b)} +\frac{b(c-a)}{a(c-b)}Q(a,b,c\mid x).\]
\end{lem}
\begin{proof}
The formula is an immediate consequence of the contiguous relation
\[c(a-b)\shyper{a,b}cx + b(c-a)\shyper {a,b+1}{c+1}x+a(b-c)\shyper{a+1,b}{c+1}x=0,\]
which is easily proved by equating coefficients of powers of $x$.
\end{proof}

Equivalently, Lemma \ref{l-Q} asserts that $ca+b(c-a)Q(a,b,c\mid x)$ is symmetric in $a$ and $b$.

\begin{prop}
With $A(x)$ as in Lemma \ref{l-hankel},
we have
\[H_n^{(2)}(A)=\left(\frac{a(c-b)}{c(a-b)}\frac{\r{a+1}n
\r{c-b+1}n}{\r{b+1}n \r{c-a+1}n}-\frac {b(c-a)}{c(a-b)}\right)H_{n+1}(A).\]
\end{prop}

\begin{proof}
First note that if $u(x)=\alpha+\beta v(x)$, where $\alpha$ and $\beta$ are constants,
then
\[H_{n+1}(u)=\beta^{n+1}H_{n+1}(v)+\alpha\beta^{n}H_{n}^{(2)}(v),\]
so
\begin{equation}
H_{n}^{(2)}(v)=\frac 1{\alpha\beta^{n}}H_{n+1}(u) -\frac\beta\alpha
H_{n+1}(v).
\label{e-H2a}
\end{equation}
Now take $u=Q(b,a,c \mid x)$ and  $v=Q(a,b,c\mid x)$,
so that  $u=\alpha + \beta v$ by Lemma \ref{l-Q},
 where $\alpha=c(a-b)/a(c-b)$ and $\beta=b(c-a)/a(c-b)$. Then by Lemma
\ref{l-hankel}, we have
\begin{align}
\frac{H_{n+1}(u)}{H_{n+1}(v)}&=\prod_{i=1}^n
  \frac{\r{a+1}i}{\r ai}
  \frac {\r bi}{\r{b+1}i}
  \frac {\r{c-b+1}i}{\r{c-b}i}
  \frac  {\r{c-a}i}{\r{c-a+1}i}\nonumber\\
  &=\prod_{i=1}^n \frac{b(c-a)}{a(c-b)}\frac{(a+i)(c-b+i)}{(b+i)(c-a+i)}
  =\left[\frac{b(c-a)}{a(c-b)}\right]^n \frac{\r {a+1}n \r{c-b+1}n}{\r {b+1}n
\r{c-a+1}n},\label{e-H2b}
\end{align}
and by \eqref{e-H2a} we have
\begin{equation}
\frac{H_{n}^{(2)}(v)}{H_{n+1}(v)}
  =\frac {a(c-b)}{c(a-b)}\left[\frac{a(c-b)}{b(c-a)}\right]^n\frac{H_{n+1}(u)}{H_{n+1}(v)} -
  \frac{b(c-a)}{c(a-b)}.\label{e-H2c}
\end{equation}
The result follows from \eqref{e-H2b} and \eqref{e-H2c}.
\end{proof}

Tamm \cite{tamm} evaluated the determinants $U_n$ and $V_n$ by first showing that
\begin{equation}\label{e-gmain}
\sum_{n=0}^\infty a_n x^n  = \thyper{\frac23,\frac43}{\frac32}{ \frac{27}{4} x  }
  \bigg/  \thyper{\frac23,\frac13}{\frac12}{\frac{27}{4} x  }.
\end{equation}
Given \eqref{e-gmain}, it follows from Lemma \ref{l-hankel} that
\begin{align*}
U_n&=\prod_{i=1}^{n-1}\frac
 {\r{\frac23}i \r{\frac16}i \r{\frac43}i \r{\frac 56}i}
  {\r{\frac12}{2i} \r{\frac32}{2i}}
\left(\tf\right)^{2i}
\\
\intertext{and}
V_n&=\prod_{i=0}^{n}
\frac23
\frac
 {\r{\frac23}i \r{\frac16}i \r{\frac13}i \r{-\frac 16}i}
  {\r{\frac12}{2i} \r{-\frac12}{2i}}
\left(\tf\right)^{2i}
\end{align*}
So \eqref{e-typea} and \eqref{e-typeb} will follow from
\begin{align}
\frac
 {\r{\frac23}i \r{\frac16}i \r{\frac43}i \r{\frac 56}i}
  {\r{\frac12}{2i} \r{\frac32}{2i}}
 \left(\tf\right)^{2i}
&=\frac{(3i+1)(6i)!\,(2i)!}{(4i+1)!\,(4i)!}\label{e-id1}\\
\intertext{and}
\frac23
\frac
 {\r{\frac23}i \r{\frac16}i \r{\frac13}i \r{-\frac 16}i}
  {\r{\frac12}{2i} \r{-\frac12}{2i}}
\left(\tf\right)^{2i}
&=
\dfrac{\displaystyle\binom {6i-2}{2i}}{\displaystyle 2\binom{4i-1}{2i}}
\label{e-id2}
\end{align}
for $i\ge1$. These identities are most easily verified by using the fact that if
$A_1=B_1$ and
$A_{i+1}/A_i=B_{i+1}/B_i$ for $i\ge1$, then $A_i=B_i$ for all $i\ge1$. It is
interesting to note that although \eqref{e-id1} holds for $i=0$, \eqref{e-id2}
does not.

\section{Hypergeometric series evaluations}

Let $f=g-1=\sum_{n=1}^\infty a_n x^n=\sum_{n=1}^\infty \frac1{2n+1}\binom{3n}{n} x^n$. In this section we
study cases of Gauss's continued fraction
\eqref{e-gausscf} that can be expressed in terms of $f$. We found empirically that there are ten cases of
\eqref{e-gausscf} that can be expressed as polynomials in $f$. We believe there are no others, but we
do not have a proof of this. Since $a\ne b$ in all of these cases, by Lemma \ref{l-Q} they
must come in pairs which are the same, except for their constant terms, up to a constant factor. It
turns out that one element of each of these pairs factors as $(1+f)(1+rf)$, where $r$ is 0, 1,
$\frac12$,
$-\frac12$, or
$\frac25$, while the other  does not factor nicely. We have no explanation for this phenomenon.

Note that \eqref{e-cf1a} is the same as \eqref{e-gmain}.

\begin{thm}
\label{t-cfs}
\begin{subequations}
 We have the following cases of Gauss's continued fraction:
\begin{align}
1+f  &= \thyper{\frac23,\frac43}{\frac32}{ \frac{27}{4} x  }
  \bigg/  \thyper{\frac23,\frac13}{\frac12}{\frac{27}{4} x  }\label{e-cf1a}\\
(1+f)^2 &=\thyper{\frac43,\frac53}{\frac52}{\frac{27}{4}x}\bigg/
\thyper{\frac43,\frac23}{\frac32}{\frac{27}{4}x}\label{e-cf2a} \\
(1+f)(1+\tfrac12f) &=\thyper{\frac53,\frac73}{\frac72}{\frac{27}{4}x}\bigg/
\thyper{\frac53,\frac43}{\frac52}{\frac{27}{4}x}\label{e-cf3a}\\
(1+f)(1-\tfrac12f)&=\thyper{\frac53,\frac73}{\frac52}{\frac{27}{4}x}\bigg/
\thyper{\frac53,\frac43}{\frac32}{\frac{27}{4}x}\label{e-cf4a} \\
(1 +f)(1+\tfrac25f)
&=\thyper{\frac23,\frac43}{\frac52}{\frac{27}{4}x} \bigg/
\thyper{\frac23,\frac13}{\frac32}{\frac{27}{4}x}\label{e-cf5a}
\end{align}
\end{subequations}
Their companions are
\begin{subequations}
\begin{align}
1-\tfrac12 f &=\thyper{\frac13,\frac53}{\frac32}{\frac{27}{4}x} \bigg/
 \thyper{\frac13,\frac23}{\frac12}{\frac{27}{4}x}\label{e-cf1b} \\
 1+\tfrac15 f +\tfrac{1}{10}f^2 &=\thyper{\frac23,\frac73}{\frac52}{\frac{27}{4}x}\bigg/
\thyper{\frac23,\frac43}{\frac32}{\frac{27}{4}x}\label{e-cf2b}\\
 1+\tfrac67f +\tfrac27 f^2 &=\thyper{\frac43,\frac83}{\frac72}{\frac{27}{4}x} \bigg/
\thyper{\frac43,\frac53}{\frac52}{\frac{27}{4}x}\label{e-cf3b}\\
1-\tfrac25 f +\tfrac25 f^2
&=\thyper{\frac43,\frac83}{\frac52}{\frac{27}{4}x}\bigg/
\thyper{\frac43,\frac53}{\frac32}{\frac{27}{4}x}\label{e-cf4b} \\
1+\tfrac12f +\tfrac17 f^2
&=\thyper{\frac13,\frac53}{\frac52}{\frac{27}{4}x}\bigg/
\thyper{\frac13,\frac23}{\frac32}{\frac{27}{4}x}\label{e-cf5b}
\end{align}
\end{subequations}
\end{thm}

In order to prove Theorem \ref{t-cfs}, we need formulas  for some rational functions of $f$
that are easily proved by Lagrange inversion.

\begin{lem}
Let $f=\sum_{n=1}^\infty \frac1{2n+1}\binom{3n}{n} x^n$. Then $f$ satisfies the functional equation
$f=x(1+f)^3$ and
\begin{align}
f^k &= \sum_{n=k}^\infty \frac kn\binom{3n}{ n-k} x^n
\label{e-f-kth}\\
(1+f)^k &= \sum_{n=0}^{\infty}\frac k{3n+k}\binom {3n+k}n x^n
\label{e-g-kth}\\
\frac{ (1+f)^{k+1} }{1-2f} &= \sum_{n= 0}^\infty \binom {3n+k}n  x^n. \label{e-kth}
\end{align}
In particular,
\begin{align}
1+f &=\thyper{\frac13,\frac23}{\frac32}{\frac{27}{4}x} \label{e-g-0}\\
\frac{1+f}{1-2f}
& = \thyper{\frac13,\frac23}{\frac12}{\frac{27}{4} x  }
\label{e-g-1} \\
\frac{(1+f)^2}{1-2f}& = \thyper{\frac43,\frac23}{\frac32}{\frac{27}{4} x }.
\label{e-g-2}
\end{align}
\end{lem}

\begin{proof}
We use the following form of the Lagrange inversion formula (see
\cite[Theorem 2.1]{gessel} or \cite[Theorem 1.2.4]{gj}): If $G(t)
$ is a formal power series, then there is a unique formal power
series $h=h(x)$ satisfying $h=xG(h)$, and
\begin{align}
[x^n]\, h^k &= \frac kn[t^{n-k}]\, G(t)^n,
\text{  for $n,k>0$,}\label{e-i} \\
[x^n]\,\frac{h^k}{1-xG^\prime (h)} &=
[t^{n-k}]\, G(t)^n, \text{  for $n,k\ge 0$}.\label{e-ii}
\end{align}

Let us define $f$ to be the unique formal power series satisfying $f=x(1+f)^3$. With $G(t)=(1+t)^3$, \eqref{e-i}
gives \eqref{e-f-kth},
and
the case $k=1$ gives that the coefficient of $x^n$  in $f$ for $n\ge1$ is $\frac1n\binom{3n}{n-1}=\frac{1}{2n+1}\binom{3n}n$.

Replacing with $f$ with $x(1+f)^3$ and $k$ with $j$ in \eqref{e-f-kth}, and dividing both sides by $x^j$, gives
\[
(1+f)^{3j}=\sum_{n=0}^\infty \frac j{n+j}\binom{3n+3j}{n}x^n.
\]
Since the coefficient of $x^n$ on each side is a polynomial in $j$, we may set $j=k/3$ to obtain \eqref{e-g-kth}.

{}From \eqref{e-ii} we have
$$\frac{f^j}{1-3x(1+f)^2} = \sum_{n=j}^\infty {3n\choose n-j} x^n.$$
Replacing $f$ by $x(1+f)^3$ in the numerator, and replacing $x(1+f)^2$ by
$f/(1+f)$ in the denominator, gives
\begin{align*}
\frac{x^j (1+f)^{3j+1}
}{1-2f} & = \sum_{n=j}^{\infty}{3n\choose n-j} x^n \\
  &=
\sum_{n=0}^\infty {3n+3j\choose n}
x^{n+j},
\end{align*}
so
$$\frac{(1+f)^{3j+1}}{1-2f} = \sum_{n=0}^\infty {3n+3j
\choose n} x^n.$$
As before, we may set
$j=k/3$ to obtain
\eqref{e-kth}.
\end{proof}

\begin{proof}[Proof of Theorem \ref{t-cfs}]
Formulas \eqref{e-cf1a}--\eqref{e-cf5a} follow from the evaluations of their numerators and denominators:
\eqref{e-g-0}, \eqref{e-g-1}, \eqref{e-g-2}, and
\begin{align}
\thyper{\frac23,\frac43}{\frac52}{\frac{27}{4}x}
 &=(1+f)^2(1+\textstyle\frac25f)\label{e-q41}\\
\thyper{\frac43,\frac53}{\frac52}{\frac{27}{4}x}
 &=\frac{(1+f)^4}{1-2f}\label{e-q42}\\
\thyper{\frac53,\frac73}{\frac72}{\frac{27}{4}x}
&=\frac{(1+f)^5(1+\frac12f)}{1-2f}\label{e-q43}\\
\thyper{\frac43,\frac53}{\frac32}{\frac{27}{4}x}
&=\frac{(1+f)^4}{(1-2f)^3}\label{e-q44}\\
\thyper{\frac53,\frac73}{\frac52}{\frac{27}{4}x}
&=\frac{(1+f)^5(1-\frac12f)}{(1-2f)^3}.\label{e-q45}
\end{align}
Our original
derivations of these formulas were through the $_2F_1$ contiguous relations
\cite[p.~558]{AS}, but once we have found
them, we can verify
\eqref{e-q41}--\eqref{e-q43} by
by taking appropriate linear
combinations of \eqref{e-f-kth} and
\eqref{e-kth}.  Formulas
 \eqref{e-q44} and \eqref{e-q45}
can be proved by applying the formula
\[\shyper{a+1,b+1}{c+1}x=\frac{c}{ab}\frac{d\ }{dx}\shyper {a,b}cx\] to \eqref{e-g-1}
and \eqref{e-g-2} and using the fact
that
$df/dx=(1+f)^4/(1-2f)$.

Formulas \eqref{e-cf1b}--\eqref{e-cf5b} can be proved similarly; alternatively, they can be derived from
 \eqref{e-cf1a}--\eqref{e-cf5a} by using Lemma \ref{l-Q}.
\end{proof}

Now we apply Lemma \ref{l-hankel} to the formulas of Theorem \ref{t-cfs}. First we normalize the
coefficient sequences that occur in \eqref{e-cf1a}--\eqref{e-cf5a} to make them integers, using
\eqref{e-f-kth} to find formulas for the coefficients. We define the sequence $a_n$, $b_n$, $c_n$,
$d_n$, and $e_n$ by
\begin{align*}
1+f&=\sum_{n=0}^\infty a_n x^n         &a_n&=\frac{(3n)!}{n!\,(2n+1)!}=\frac1{2n+1}\binom{3n}n\\
(1+f)^2&=\sum_{n=0}^\infty b_n x^n     &b_n&=\frac{(3n+1)!}{(n+1)!\,(2n+1)!}=\frac1{n+1}\binom{3n+1}n\\
(1+f)(2+f)&=\sum_{n=0}^\infty c_n x^n  &c_n&=2\frac{(3n)!}{(n+1)!\,(2n)!}=\frac2{n+1}\binom{3n}n=a_n+b_n\\
(1+f)(2-f)&=\sum_{n=0}^\infty d_n x^n  &d_n&=2\frac{(3n)!}{(n+1)!\,(2n+1)!}=3a_n-b_n\\
(1+f)(5+2f)&=\sum_{n=0}^\infty e_n x^n &e_n&=(9n+5)\frac{(3n)!}{(n+1)!\,(2n+1)!}=3a_n+2b_n
\end{align*}

Here is a table of the first few values of these numbers

$$\begin{array}{r|rrrrrrrr}
n&0\hfil&1\hfil&2&3\hfil&\hfil4\hfil&5\hfil&6\hfil&7\hfil\\  \hline
a_n&1&1&3&12&55&273&1428&7752\\
b_n&1&2&7&30&143&728&3876&21318\\
c_n&2&3&10&42&198&1001&5304&29070\\
d_n&2&1&2&6&22&91&408&1938\\
e_n&5&7&23&96&451&2275&12036&65892
\end{array}$$

The sequences $a_n$ and $b_n$ are well-known, and have simple
combinatorial interpretations in terms of lattice paths: $a_n$ is
the number of paths, with steps $(1,0)$ and $(0,1)$, from $(0,0)$
to $(2n,n)$ that never rise above (but may touch) the line $x=2y$
and $b_n$ is the number of paths from $(0,0)$ to $(2n,n)$ that
never rise above (but may touch) the line $x=2y-1$ (see, e.g.,
Gessel \cite{ira-path}).  Moreover, for $n>0$, $d_n$ is the number
of two-stack-sortable permutations of $\{1,2,\ldots, n\}$. (See,
e.g., West \cite{west} and Zeilberger \cite{zeil}.) The sequences
$c_n$ and $e_n$ are apparently not well-known.

Let us write $H_n(a)$ for $H_n\bigl(\sum_{n=0}^\infty a_n x^n\bigr)$, and similarly for other letters replacing $a$. Then
applying  Lemma \ref{l-hankel} and Theorem \ref{t-cfs}gives
\begin{align*}
H_n(a)&=
\prod_{i=0}^{n-1}
\frac{(\frac{2}{3})_{i}(\frac{1}{6})_{i}(\frac{4}{3})_{i}(\frac{5}{6})_{i}}
{(\frac{1}{2})_{2\,i}(\frac{3}{2})_{2\,i}}
\left(\tf\right)^{2i}\\
H^1_n(a)&=
\prod_{i=1}^{n}\frac{2}{3}
\frac{(\frac{2}{3})_{i}(\frac{1}{6})_{i}(\frac{1}{3})_{i}(-\frac{1}{6})_{i}}
{(\frac{1}{2})_{2\,i}(-\frac{1}{2})_{2\,i}}
\left(\tf\right)^{2i}\\
H_n(b)&=\prod_{i=0}^{n-1}\frac{(\frac{4}{3})_{i}(\frac{5}{6})_{i}(\frac{5}{3})_{i}(\frac{7}{6})_{i}}
{(\frac{3}{2})_{2\,i}(\frac{5}{2})_{2\,i}}
\left(\tf\right)^{2i}\\
H_n^1(b)&=
\prod_{i=1}^{n}
\frac{(\frac{4}{3})_{i}(\frac{5}{6})_{i}(\frac{2}{3})_{i}(\frac{1}{6})_{i}}
{(\frac{3}{2})_{2\,i}(\frac{1}{2})_{2\,i}}
\left(\tf\right)^{2i}\\
H_n(c)&=
\prod_{i=0}^{n-1}
2\frac{(\frac{5}{3})_{i}(\frac{7}{6})_{i}(\frac{7}{3})_{i}(\frac{11}{6})_{i}}
{(\frac{5}{2})_{2\,i}(\frac{7}{2})_{2\,i}}
\left(\tf\right)^{2i}\\
H_n^1(c)&=
\prod_{i=1}^{n}
\frac{(\frac{5}{3})_{i}(\frac{7}{6})_{i}(\frac{4}{3})_{i}(\frac{5}{6})_{i}}
{(\frac{5}{2})_{2\,i}(\frac{3}{2})_{2\,i}}
\left(\tf\right)^{2i}\\
H_n(d)&=
\prod_{i=0}^{n-1}
2\frac{(\frac{5}{3})_{i}(\frac{1}{6})_{i}(\frac{7}{3})_{i}(\frac{5}{6})_{i}}
{(\frac{5}{2})_{2\,i}(\frac{3}{2})_{2\,i}}
\left(\tf\right)^{2i}\\
H_n^1(d)&=
(-1)^n
\prod_{i=1}^{n}
\frac{(\frac{5}{3})_{i}(\frac{1}{6})_{i}(\frac{4}{3})_{i}(-\frac{1}{6})_{i}}
{(\frac{3}{2})_{2\,i}(\frac{1}{2})_{2\,i}}
\left(\tf\right)^{2i}\\
H_n(e)&=
\prod_{i=0}^{n-1} 5
\frac{(\frac{2}{3})_{i}(\frac{7}{6})_{i}(\frac{4}{3})_{i}(\frac{11}{6})_{i}}
{(\frac{3}{2})_{2\,i}(\frac{5}{2})_{2\,i}}
\left(\tf\right)^{2i}\\
H_n^1(e)&=
\prod_{i=1}^{n}2
\frac{(\frac{2}{3})_{i}(\frac{7}{6})_{i}(\frac{1}{3})_{i}(\frac{5}{6})_{i}}
{(\frac{3}{2})_{2\,i}(\frac{1}{2})_{2\,i}}
\left(\tf\right)^{2i}
\end{align*}

Here is a table of the values of these Hankel determinants:

$$\begin{array}{r|rrrrrrr}
n&1\hfil&2&3\ &4\ &5\ \ &6\ \ \ \ &7\ \ \ \ \ \ \ \\  \hline
H_n(a)&1&2&11&170&7429&920460&323801820\\
H_n^1(a)&1&3&26&646&45885&9304650&5382618660\\
H_n(b)&1&3&26&646&45885&9304650&5382618660\\
H_n^1(b)&2&11&170&7429&920460&323801820&323674802088\\
H_n(c)&2&11&170&7429&920460&323801820&323674802088\\
H_n^1(c)&3&26&646&45885&9304650&5382618660&8878734657276\\
H_n(d)&2&3&10&85&1932&120060&20648232\\
H_n^1(d)&1&2&10&133&4830&485460&136112196\\
H_n(e)&5&66&2431&252586&74327145&62062015500&147198472495020\\
H_n^1(e)&7&143&8398&1411510&677688675&928501718850&3628173844041420\\
\end{array}$$

It is apparent from the table that
\begin{align}
\label{e-U-abc} U_n=H_n(a)=H_{n-1}^1(b)=H_{n-1}(c),
\end{align} and that
$V_n=H_n^1(a)=H_n(b)=H_{n-1}^1(c)$, and these are easily verified
from the formulas. The combinatorial interpretations of $U_n$ and
$V_n$ have already been discussed. The numbers $H_n(e)$ were shown
by Kuperberg \cite[Theorem 5]{kuperberg} to count certain
alternating sign matrices. In Kuperberg's notation,
$H_n(e)=A_{\text{UU}}^{(2)}(4n;1,1,1)$.

There are also Hankel determinant evaluations corresponding to
\eqref{e-cf1b}--\eqref{e-cf5b}, normalized to make the entries
integers. These evaluations can be found in Krattenthaler \cite[Theorem 30]{kratt2}.

\section{Determinants and Two-Variable
Generating Functions}

In this section we describe a method for transforming determinants whose
entries are given as coefficients of generating functions. (A related approach was used in \cite{gesselprob} to evaluate Hankel determinants of Bell numbers.) Using
this technique, we are able to convert the determinants for $U_n$ and
$V_n$ in \eqref{e-s1} and \eqref{e-s2}  into the known determinant evaluations given in \eqref{e-un} and
\eqref{e-vn}. (Conversely, the evaluations of these Hankel determinants give new proofs of \eqref{e-un} and
\eqref{e-vn}.)
These two determinants are
special cases of a determinant evaluation of Mills, Robbins, and Rumsey
\cite{tcsymm} (see \cite[Theorem~37]{kratt}
for related determinants):
\begin{multline}
\detijn{ \binom{i+j+r}{2i - j}}\\=(-1)^{\chi (n\equiv 3 \bmod
4)}2^{\binom{n-1}{2}} \prod_{i=1}^{n-1}\frac{(r+i+1)_{\lfloor
(i+1)/2 \rfloor}(-r-3n+i+\frac{3}{2})_{\lfloor i/2
\rfloor}}{(i)_i},\label{e-gen-det}
\end{multline}
where $\chi(S)=1$ if $S$ is true and $\chi(S)=0$ otherwise. There
exist short direct proofs of \eqref{e-gen-det} (see
\cite{AB-gen-det, K-gen-det, PW-gen-det}), but no really simple
proof.

Suppose that we have a two-variable generating
function
\[ D(x,y) = \sum_{i,j=0}^\infty d_{i,j} x^i y^j. \]
Let $[D(x,y)]_n$ be
the determinant of the $n \times n$ matrix
\[    (d_{i,j})_{0 \leq i,j \leq
n-1} .\]

The following  rules can be used to transform the
determinant $[D(x,y)]_n$ to a determinant with the same value:

\smallskip
\noindent\textbf {Constant Rules.}
Let $c$ be a non-zero constant. Then
\begin{align*}
[cD(x,y)]_n &= c^n [D(x,y)]_n, \\ 
\intertext{and}
[D(cx,y)]_n &= c^{n\choose 2} [D(x,y)]_n. 
\end{align*}
\textbf {Product Rule.}
If $u(x)$ is any formal power series with $u(0)=1$,
then
\[[u(x)D(x,y)]_n=[D(x,y)]_n.\]
\textbf {Composition Rule.}
If $v(x)$ is any formal power series with $v(0)=0$
and
$v'(0)=1$, then
$$[D(v(x),y)]_n =[D(x,y)]_n.$$

\smallskip
The product and composition rules hold because the transformed
determinants are obtained from the original determinants by
elementary row operations. Equivalently, the new matrix is
obtained by multiplying the old matrix on the left by a
matrix with determinant
$1$. Note that all of these transformations can be applied to $y$
as well as to $x$.

The Hankel determinants $H_n(A)$ and $H_n^1(A)$ of a formal power
series $A(x)$ are given by
\begin{align}
H_n(A) &=\left[\frac{xA(x) -
yA(y)}{x-y}\right]_n, \label{e-hankel0-gen}\\
H_n^1(A) &= \left[\frac{A(x) -
A(y)}{x-y}\right]_n.\label{e-hankel1-gen}
\end{align}

\begin{proof}[Proof of \eqref{e-un} and \eqref{e-vn}] The generating
function for the Hankel determinant $H_n(g)$ is
\begin{equation}
\label{e-ghankel}
\frac{xg(x) - yg(y)}{x-y}.
\end{equation}
Since $f/(1+f)^3=x$, $f$ is the compositional inverse of $x/(1+x)^3$, and thus $f\bigl(x/(1+x)^3\bigr) = x$. Since $g=1+f$, we have
$g\bigl(x/(1+x)^3\bigr)=1+x$.

Now let us substitute
$x\to x/(1+x)^3$, $y \to y/(1+y)^3$ in \eqref{e-ghankel}.
After simplifying, we obtain
\[\frac{(1 - xy)(1 + x) (1 + y)}{1 - x y^2 - 3xy - x^2 y}.\]
Then dividing by $(1+x)(1+y)$, we get
$$\frac{1 - xy}{1 - xy^2 - 3xy -
x^2y}.$$
Next, we show that
\begin{equation}
\frac{1 - xy}{1 - xy^2 - 3xy -
x^2y}= \sum_{i,j} {i + j \choose 2i - j} x^i y^j.
\label{e-gf1}
\end{equation}
Multiplying both sides of \eqref{e-gf1} by
${1 -
xy^2 - 3xy - x^2y}$ and equating coefficients of $x^my^n$ shows that
\eqref{e-gf1}
is equivalent to the recurrence
\[{m+n\choose 2m-n}-{m+n-3\choose
2m-n}-3{m+n-2\choose 2m-n-1}
-{m+n-3\choose 2m-n-3}=
   \begin{cases}

\phantom{-}1,&\text{if $m=n=0$}\\
   -1,&\text{if $m=n=1$}\\

\phantom{-}0,&\text{otherwise}
   \end{cases}
\]
where we interpret the binomial coefficient $a\choose b$ as 0 if
either $a$ or $b$ is negative, and the verification of the
recurrence is straightforward. (We will give another proof of
\eqref{e-gf1} in Example \ref{ex-1}.)
This completes the proof of
\eqref{e-un}.

For equation \eqref{e-vn}, we need to consider the generating
function
$$(g(x)
- g(y))/(x-y).$$ Making the same substitution as before gives
$$ \frac{(1+x)^3
(1+y)^3}{1 - xy^2 - 3xy - x^2y}.$$
Dividing by $(1+x)^2(1+y)^3$
gives
$$\frac{1+x}{1 - x y^2 - 3xy - x^2 y},$$
which can be shown, by the same
method as before, to equal
\[\sum_{i,j} {i + j + 1 \choose 2i - j} x^i y^j.
\qedhere
\]
\end{proof}

To transform in this way the more general determinant on the left side of \eqref{e-gen-det}, we would start with the generating
function
\begin{align}\label{e-r}
\sum_{i,j}{i+j+r \choose 2i -
j}x^iy^j=
\frac{\displaystyle\sum_{n=0}^\infty \left[{r+n\choose
2n}+{r+n-2\choose 2n-1}y\right]x^n}
{1 - x y^2 - 3xy - x^2 y}.
\end{align}
The generating function in $r$ of \eqref{e-r} is derived in
\eqref{e-gfzzz}. The sums in the numerator can be evaluated
explicitly by making an appropriate substitution in the identities
\begin{align*}\sum_{n=0}^\infty{r+n\choose
2n} (-4\sin^2\theta)^n&={\cos(2r+1)\theta\over
\cos\theta}\\
\sum_{n=0}^\infty
\binom{r+n-2}{2n-1}(-4\sin^2\theta)^n&=-2\tan\theta\,\sin
2(r-1)\theta.
\end{align*}
However we have not been able to use these formulas to prove
\eqref{e-gen-det}.

Another application of this method gives a family of generating functions that have the same Hankel determinants.

\begin{thm}
\label{t-same}
Let $A(x)$ be a formal power series with $A(0)=1$ and let $c$ be a
constant. Then we have
\begin{align}
H_n\left(\frac{A(x)}{1-cxA(x)}\right)&=H_n(A)\label{e-hna1}
\end{align}
for all $n$, and
\begin{align}\label{e-hna11}
H_n\left(\frac{1}{1-cxA(x)}\right) &= c^{n-1}H_{n-1}^1(A)
\end{align}
for  $n\ge 1$.
\end{thm}
\begin{proof}
We use the method of generating functions to evaluate these
determinants. By \eqref{e-hankel0-gen},
\begin{align*}
H_n\left(\frac{A(x)}{1-cxA(x)}\right)
&=\begin{bmatrix}\frac{\displaystyle\frac{xA(x)}{1-cxA(x)}-
   \frac{yA(y)}{1-cyA(y)}}{\strut\displaystyle x-y}\end{bmatrix}_n \\
&=\left[\frac{1}{(1-cxA(x))(1-cyA(y))}
\frac{xA(x)-yA(y)}{x-y}\right]_n.
\end{align*}
Since$(1-cxA(x))^{-1}$ is a formal power series with constant term
$1$, we get
\begin{align*}
H_n\left(\frac{A(x)}{1-cxA(x)}\right)
&=\left[\frac{xA(x)-yA(y)}{x-y}\right]_n =H_n(A).
\end{align*}

A similar computation shows that
\begin{align*}
H_n\left(\frac{1}{1-cxA(x)}\right) &=\left[1+cxy
\frac{A(x)-A(y)}{x-y}\right]_n \\
&=\left[c \frac{A(x)-A(y)}{x-y}\right]_{n-1}=c^{n-1}H_{n-1}^1(A),
\end{align*}
since $\left[1+xy D(x,y)\right]_n$ is the determinant of a block matrix of two
blocks, with the first block $[1]$ and the second block
$\left[D(x,y)\right]_{n-1}$. \end{proof}

We now prove \eqref{e-err} and \eqref{e-sdet}.
First we set $c=u-1$ and $A=f/x$ in \eqref{e-hna1}, getting
$$V_n=\detijn{a_{i+j+1}}=H_n(f/x)=H_n\left(\frac{f/x}{1+(1-u)f}\right).$$
Next we show that
\begin{equation}
\label{e-fs}
\frac{f/x}{1+(1-u)f}=\sum_{n=0}^\infty s_n(u)x^n,
\end{equation}
where
\begin{equation}
\label{e-su}
s_n(u)=\sum_{k=0}^{n}\frac{k+1}{n+1}\binom{3n-k+1}{n-k}u^k.
\end{equation}
We have
\begin{align*}
\frac{f/x}{1+(1-u)f} &= \frac{f/x}{1+f}\cdot\frac 1 {1-uf/(1+f)}\\
  &=\frac {(1+f)^2}{1-ux(1+f)^2}, \text{ since $f=x(1+f)^3$,}\\
  &=\sum_{k=0}^\infty u^k x^k (1+f)^{2k+2}\\
  &=\sum_{k=0}^\infty u^k x^k \sum_{m=0}^\infty \frac{2k+2}{3m+2k+2}
   \binom{3m+2k+2}m x^m, \text{ by \eqref{e-g-kth},}\\
  &=\sum_{n=0}^\infty x^n\sum_{k=0}^n\frac{k+1}{n+1}\binom{3n-k+1}{n-k}u^k,
\end{align*}
which proves \eqref{e-fs}. Then $s_n(1)= a_{n+1}$ from \eqref{e-fs}, $s_n(0)=\frac1{n+1}\binom{3n+1}n$ by setting $u=0$ in \eqref{e-su}, and
$s_n(3) = \binom{3n+2}{n}$ follows from \eqref{e-kth}. This completes the proof of
\eqref{e-err}.

Next we prove \eqref{e-sdet}, which by \eqref{e-fs} is equivalent to
\begin{equation}
\label{e-uf}
H_n\left(u^{-1}+\frac{f}{1+(1-u)f} \right) =U_n/u.
\end{equation}
We have
$$ u^{-1}+\frac{f}{1+(1-u)f}=  \frac{u^{-1}}{1-ux(1+f)^2},  $$
so by \eqref{e-hna11}, the Hankel determinant is equal to $u^{-1} H_{n-1}^1((1+f)^2)$. In the notation of Section 3, this is
$u^{-1} H_{n-1}^1(b)$, which by  \eqref{e-U-abc} is equal to
$u^{-1} U_n$.
\smallskip

We also have an analogue of Theorem \ref{t-same} for the Hankel determinants $H_n^1$.
\begin{thm}
Let $A(x)$ be a formal power series with $A(0)=1$ and let $c\neq 1$ be a constant.
Then we have
\begin{align}
H_n^1\left(\frac{A(x)}{1-cA(x)}\right)&=
 (1-c)^{-2n} H_n^1(A)
\label{e-hn1a}
\end{align}
\end{thm}
\begin{proof}
We use the method of generating functions. By
\eqref{e-hankel1-gen},
\begin{align*}
H_n^1\left(\frac{A(x)}{1-cA(x)}\right)
&=\begin{bmatrix}\frac{\displaystyle\frac{A(x)}{1-cA(x)}-
\frac{A(y)}{1-cA(y)}}{\displaystyle \strut x-y}\end{bmatrix}_n \\
&=\left[\frac{1}{(1-cA(x))(1-cA(y))}
\frac{A(x)-A(y)}{x-y}\right]_n.
\end{align*}
Since $(1-cA(x))^{-1}$ is a formal power series with
constant term $(1-c)^{-1}$ when $c\neq 1$, we get
\begin{align*}
H_n^1\left(\frac{(1-c)^2A(x)}{1-cA(x)}\right)
&=\left[(1-c)^{-2}\frac{A(x)-A(y)}{x-y}\right]_n =(1-c)^{-2n}H_n^1(A).
\qedhere
\end{align*}
\end{proof}

\section{A Hankel Determinant for the Number of Alternating Sign Matrices}

Let $\mathcal{A}_n$ be the number of $n\times n$ alternating sign
matrices. It is well-known that
$$\mathcal{A}_n =
\prod_{k=0}^{n-1} \frac{(3k+1)!}{(n+k)!},$$
as conjectured by Mills Robbins and Rumsey \cite{mrr} and proved by Zeilberger \cite{zeil-alt} and
Kuperberg \cite{kuperberg-alt}.

The numbers $\mathcal{A}_n$ also count totally
symmetric, self-complementary plane partitions, as shown
by Andrews \cite{ge-andrew}. We find, up to a power of
$3$, a Hankel determinant expression for $\mathcal{A}_n$.

Let
\begin{align}
\C(x)& =
\frac{1-(1-9 x)^{{1/ 3}} } {3x}.
\end{align}
The coefficients of $\C(x)$ are positive integers that are analogous to Catalan numbers. They have no known combinatorial interpretation and have been little studied, but they do appear in  \cite[Eq.
61]{lang}.
\begin{thm}
\label{t-asm}
The number of $n\times n$ alternating sign matrices is
\begin{align}
\mathcal{A}_n &= 3^{-{n\choose 2}} H_n (\C).
\end{align}
\end{thm}
\begin{proof}
Let
$$D(x,y)=(x\C(x)-y\C(y))/(x-y)$$
be
the generating function for the Hankel determinant
$H_n(\C)$. It is easy to see that
$D(x/\sqrt{3},y/\sqrt{3})$ is the generating function for
$3^{-{n\choose 2}} H_n
(\C)$. We make the substitution $x\to x - \sqrt{3} x^{2}+ x^{3}$,  $y\to y - \sqrt{3}
y^{2}+ y^{3}$ in $D(x/\sqrt{3},y/\sqrt{3})$,
and simplify. The generating function becomes
$$
\frac{1}{1-\sqrt{3}(x+y)+x^2+xy+y^2}.$$
Let $\omega =-\frac12-{\sqrt{-3}\over
2}$ be a cube root of unity.
Make another substitution $x\to
-\sqrt{-1}x/(1+\omega x), y\to \sqrt{-1}y/(1+ \omega^2 y)$, and simplify.
The
generating function becomes
$$\frac{(1+\omega x)^2\, (1+ \omega^2
y)^2}{(1-xy)(1-x-y)}.$$
Dividing by $(1+\omega x)^2\, (1+ \omega^2 y)^2$, the
generating
function becomes
$$\frac{1}{(1-x-y)(1-xy)}.$$
Multiplying by
$(1-x+x^2)(1-y)/(1-x)$, we get
$$\frac{(1-x+x^2)(1-y)}{(1-x)(1-x-y)(1-xy)} =
{x\over y(1-x-y) }+ {1\over 1-xy} -{x\over y(1-x)}.$$ Expanding
the right-hand side of the above equation, we get
$$ H_n(\C) =3^{\binom{n}{2}} \det \left({i+j\choose
i-1}+\delta_{i,j}\right)_{0\leq i,j\leq n-1},$$ where
$\delta_{i,j}$ equals $1$ if $i=j$ and $0$ otherwise. The theorem
then follows from a known formula for $\mathcal{A}_n$
\cite[p.~22]{paproof}.
\end{proof}

\begin{rem} We have another determinant expression
$$\mathcal{A}_n =
\det \left({i+j\choose i}-\delta_{i,j+1}\right),$$
since
$$\frac{1}{(1-x-y)(1-xy)}=\frac{1}{1-y+y^2} \left({1\over 1-x-y}-{y\over
1-xy}\right).$$
\end{rem}

There is a result similar to Theorem \ref{t-asm}
\begin{align*}
\label{e-C1}
\C_1 (x)& = \frac{1-(1-9 x)^{2/3} }
{3x}.
\end{align*}
Let $\mathcal{A}'_n$ be the number of cyclically symmetric plane
partitions in the $n$-cube. We have
\begin{thm}
\begin{align}
\mathcal{A}'_n &= 3^{-{n\choose 2}} H_n (\C_1).
\end{align}
\end{thm}
\begin{proof}
Let
$$D(x,y)=(x\C_1(x)-y\C_1(y))/(x-y)$$
be the generating function for the Hankel
determinant
$H_n(\C_1)$. Similarly $D(x/\sqrt{3},y/\sqrt{3})$ is the generating function
for
$3^{-{n\choose 2}} H_n
(\C_1)$.
We make the same substitution (as for $H_n(\C)$) $x\to
x - \sqrt{3} x^{2}+ x^{3}$,  $y\to y - \sqrt{3} y^{2}+ y^{3}$,
and simplify. The
generating function becomes
$$
\frac{2-\sqrt{3}(x+y)}{1-\sqrt{3}(x+y)+x^2+xy+y^2}.$$
Similarly, we make another
substitution $x\to -\sqrt{-1}x/(1+\omega x), y\to \sqrt{-1}y/(1+ \omega^2 y)$,
and simplify.
The generating function becomes
$$\frac{(2-x-y-xy)(1+\omega x)\,
(1+ \omega^2 y)}{(1-xy)(1-x-y)}.$$
Dividing by $(1+\omega x)\, (1+ \omega^2 y)$,
the generating
function
becomes
$$\frac{(2-x-y-xy)}{(1-x-y)(1-xy)}=\frac{1}{1-x-y}+\frac{1}{1-xy}.$$
So
we have
$$ H_n(\C_1) = 3^{\binom{n}{2}}\det \left({i+j\choose i}+\delta_{i,j}\right)_{0\leq
i,j\leq n-1},$$ which is equal to $3^{\binom{n}{2}}
\mathcal{A}'_n$. (See \cite[p. 177, (5.28)]{paproof}.)
\end{proof}

Since $\C(x) = \shyper{\tfrac 23,1}{2}{9x}$, we can find a continued fraction for $\C(x)$ by setting $a=\frac23, b=0, c=1$ in Lemma \ref{l-gauss}, and thus evaluate the Hankel determinant for $\C(x)$ by Lemma \ref{l-Hcf}.
Similarly, since $\C_1(x) = 2\shyper{\tfrac 13,1}{2}{9x}$, we can evaluate the Hankel determinant for $\C_1(x)$ be taking
$a=\frac13, b=0, c=1$ in Lemma \ref{l-gauss}.

The Hankel determinants $H_n(\C)$ and
$H_n(\C_1)$, can also be evaluated by a
more general result (see, e.g., \cite[Theorem 26, Eq. (3.12)]{kratt}):
\begin{align}
\det \left( \binom{A}{L_i+j} \right)_{1\le i,j \le n}
= \frac{\prod_{1\le i<j\le n} (L_i-L_j) }{\prod_{1\le i\le n}
(L_i+n)! } \frac{\prod_{1\le i\le n}
(L_i+A+1)!}{\prod_{1\le i\le n}
(A+1-i)!},
\end{align}
where $L_1,\dots , L_n$ and $A$ are indeterminates, and the
factorials are interpreted using gamma functions when necessary.

Thus these calculations  give a simple method of evaluating the determinants
$$\det\left({i+j\choose
i-1}+\delta_{i,j}\right)_{0\leq i,j\leq n-1}
\text{and}\qquad
\det\left({i+j\choose
i}+\delta_{i,j}\right)_{0\leq i,j\leq n-1}.
$$
For more information on similar determinants, see
Krattenthaler \cite[Theorems 32--35]{kratt}
\cite[Section 5.5]{kratt2}.

\section{A Combinatorial Proof of \eqref{e-un}}

For the reader's convenience, we restate equation \eqref{e-un} as
follows:
\begin{align}\label{e-un-res}
\det\left(a_{i+j} \right)\ijn =\det \left(\binom{i+j}{2i-j}
\right)\ijn .
\end{align}

\begin{figure}[ht]
\begin{center}
\input{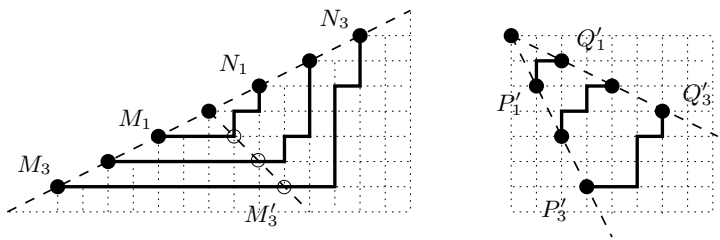}
\end{center}
 \caption{Lattice path interpretation of
\eqref{e-un-res}.\label{fig-n1}}
\end{figure}

Both sides of \eqref{e-un-res} have combinatorial meanings in
terms of nonintersecting paths (see Gessel and Viennot
\cite{gessel-viennot}). The right-hand side counts
$\mathcal{U}_R(n)$, the set of $n$-tuples of nonintersecting paths
from $P'_0,\dots ,P'_{n-1}$ to $Q'_0,\dots, Q'_{n-1}$, where
$P'_i=(i,-2i)$ and $Q'_i=(2i,-i)$. For the paths to be
nonintersecting, $P'_i$ must go to $Q'_i$. See the right picture
of Figure \ref{fig-n1}. Mills, Robbins, and Rumsey [15] in fact
gave a bijection from the type $(a)$ objects of Section 1 to such
$n$-tuples of lattice paths.

For the left-hand side, we notice that $a_n$ counts the number of
paths from $(0,0)$ to $(2n,n)$ that never go above the line
$y=x/2$. See, e.g., \cite{ira-path}. It is easy to see that the
left-hand side of \eqref{e-un-res} counts $\mathcal{U}_L(n)$, the
set of $n$-tuples of nonintersecting paths that stay below the
line $y=x/2$, from $M_0,\dots ,M_{n-1}$ to $N_0,\dots , N_{n-1}$,
where $M_i=(-2i,-i)$ and $N_i=(2i,i)$. For the paths to be
nonintersecting,
 $M_i$ must go to
$N_i$. Moreover, from the left picture of Figure \ref{fig-n1}, we
see that $M_i$ can be replaced with $M_i'=(i,-i)$.

An interesting problem is to find a bijection from
$\mathcal{U}_L(n)$ to $\mathcal{U}_R(n)$. Such a bijection will
result in a combinatorial enumeration of the type $(a)$ objects.

Both $\mathcal{U}_L(n)$ and $\mathcal{U}_R(n)$ can be easily
converted into variations of plane partitions. But we have not
found them helpful.

We find an alternative bijective proof of \eqref{e-un-res}. The
algebraic idea behind the proof is the following matrix identity
that implies \eqref{e-un-res}:
\begin{multline}
(a_{i+j})_{0\leq i,j\leq n-1} =
\left(\frac{3j+1}{3i+1}{3i+1\choose i-j}\right)_{0\leq i,j\leq
n-1}\\ \left({i+j\choose 2i-j}\right)_{0\leq i,j\leq n-1}
\left(\frac{3i+1}{3j+1}{3j+1\choose j-i}\right)_{0\leq i,j\leq
n-1},\label{e-matrix-identity}
\end{multline}
 where $$\frac{3i+1}{3j+1}{3j+1\choose j-i}=[x^i] gf^j$$
(See \eqref{e-f-kth}). Note that the left (right) transformation
matrix is a lower (upper) triangular matrix with diagonal entries
$1$. The matrix identity is obtained by carefully analyzing the
transformation we performed in Section 4 when proving
\eqref{e-un}.

The bijective proof relies on a new interpretation of $a_n$ in
terms of certain paths that we call $K$-paths. The matrix identity
\eqref{e-matrix-identity} follows easily  from the new
interpretation. This gives a bijection from $\mathcal{U}_R(n) $ to
$\mathcal{U}_{K}(n)$, the set of $n$-tuples of nonintersecting
$K$-paths resulting from the new interpretation. The desired
bijection could be completed by giving the bijection from
$\mathcal{U}_{K}(n)$ to $\mathcal{U}_L(n)$. But we have not
succeeded in this.

The new interpretation of $a_n$ consists of three kinds of paths:
normal paths, $H_2$-paths, and $V_2$-paths. A normal path has
steps $(0,1)$ and $(1,0)$. A
 path is an $H_2$ path if each horizontal step is $(2,0)$
instead of $(1,0)$. By dividing each horizontal $2$-step into two
horizontal $1$-steps, we can represent an $H_2$ path as a normal
path. Similarly, a path is a $V_2$ path if each vertical step is
$(0,2)$.

By reflecting in the line $y=-x$, we can convert an $H_2$ path
into a $V_2$ path, or a $V_2$ path into an $H_2$ path. This
bijection can convert any property of $H_2$-paths into a similar
property of $V_2$-paths.

It is well-known that the number of paths that start at $(0,0)$,
end at $(n,2n)$, and never go above the line $y=2x$ is
$a_n=\frac{1}{2n+1}{3n\choose n}$. Replacing each horizontal step
by two horizontal steps, it follows that:
\begin{prop}\label{c-h2path}
The number of $H_2$-paths $($or $V_2$-paths$)$ that start at
$(0,0)$, end at $(2n,2n)$ and never  go above the diagonal equals
$a_n$.
\end{prop}

\begin{dfn}
We call a path $P$ a $K$-path if it satisfies the following four
conditions.
\begin{enumerate}
\item[1.] The path $P$ never goes above the diagonal.
\item[2.] The part of $P$ that is below the line $y=-2x$ is a $V_2$ path.
\item[3.] The part of $P$ between the two lines $y=-2x$ and $x=-2y$
is a normal path.
\item[4.] The part of $P$ that is above the line $x=-2y$ is an $H_2$ path.
\end{enumerate}
\end{dfn}

\begin{figure}[ht]
\begin{center}
\input{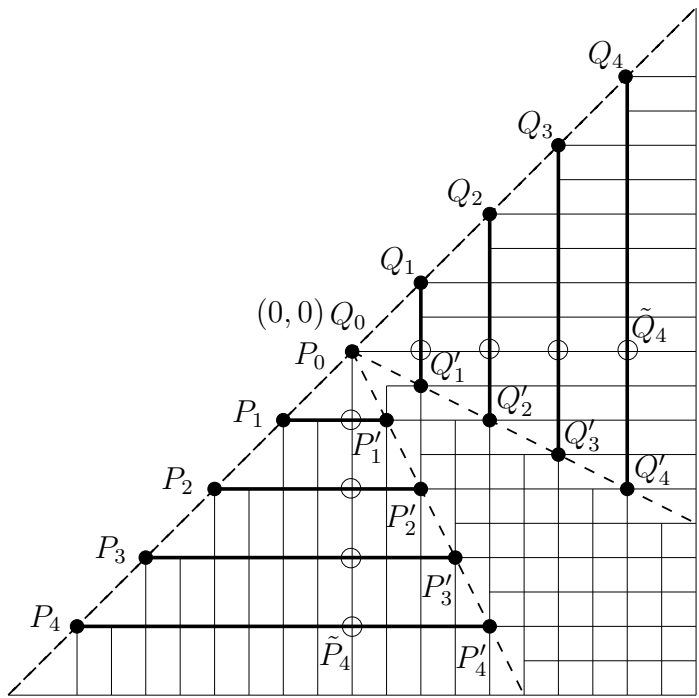}
\end{center}
\caption{The grid for $K$-paths. \label{fig-1}}
\end{figure}

{}From the definition, we  see  that a $K$-path can be uniquely
decomposed into three kinds  of paths: a $V_2$ path, followed by a
normal path, followed by an $H_2$ path. Depending on its starting
point, some of the paths may be empty. The  normal path region is
between the two lines $y=-2x$ and $x=-2y$. The steps occurring in
a $K$-path are shown in Figure \ref{fig-1}. We have

\begin{thm}\label{t-newg}
The number of $K$-paths from $(-2m,-2m)$ to $(2n,2n)$, where
$m+n\ge 0$, is $a_{m+n}$.
\end{thm}
The proof of the theorem will be given later. {}From the new
interpretation of $a_n$, $U_n$ counts $\mathcal{U}_{K}(n)$, the
set of $n$-tuples of nonintersecting $K$-paths from $P_0,\ldots
,P_{n-1}$ to $Q_0,\ldots , Q_{n-1}$, where $P_i=(-2i,-2i)$, and
$Q_i=(2i,2i)$ for $i=0,1,\ldots ,n-1$. See Figure \ref{fig-1}. For
the paths to be nonintersecting, $P_i$ must go to $Q_i$. In such
an $n$-tuples of nonintersecting $K$-paths, the path from $P_i$ to
$Q_i$ must start with a path from $P_i=(-2i,-2i)$ to
$P'_i=(i,-2i)$, and end with a path from $Q'_i=(2i,-i)$ to
$Q_i=(2i,2i)$.  So $\mathcal{U}_K(n)$ is in natural bijection with
$\mathcal{U}_R(n)$. If we count the number of $K$-paths according
to their intersections with the lines $y=-2x$ and $x=-2y$, we get
the matrix identity \eqref{e-matrix-identity}.

If $P$ is a $K$-path from $(-2m,-2m)$ to $(2n,2n)$ with $m\le 0$
(or $n\le 0$), then $P$ is an $H_2$ (or a $V_2$)-path, and Theorem
\ref{t-newg} follows from Proposition \ref{c-h2path}. So we can
assume that $m$ and $n$ are both  positive integers.

The idea of the proof of Theorem \ref{t-newg} is to show that the
number of $K$-paths from $(-2m,-2m)$ to $(2n,2n)$  is unchanged
after sliding their starting and ending points along the diagonal
by $(2,2)$.

\begin{figure}[ht]
\begin{center}
\input{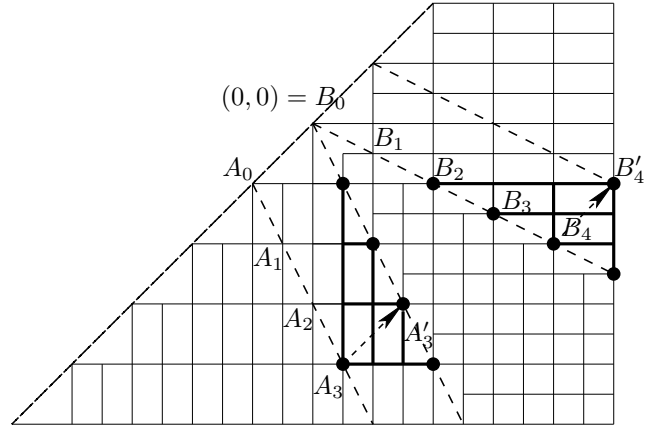}
\end{center}
\caption{Picture for the sliding lemma. \label{fig-2}}
\end{figure}

In fact, the following refinement is true. See Figure \ref{fig-2}.
\begin{lem}[Sliding Lemma]\label{l-nij}
The number of $K$-paths from $(i-2,-2i-2)$ to $(2j,-j)$ equals the
number of $K$-paths from $(i,-2i)$ to $(2j+2,-j+2)$.
\end{lem}

\begin{proof}
Let $N(i,j)$ be the number of $K$-paths from $A_i=(i-2,-2i-2)$ to
$B_j=(2j,-j)$. It is clear that $N(i,j)=0$ if $i<0$ or $j<0$.

By reflecting in the line $y=-x$, we can give a bijective proof of
the following statement: The number of $K$-paths from $(i,-2i)$ to
$(2j+2,-j+2)$ equals the number of $K$-paths from $(j-2,-2j-2)$ to
$(2i,-i)$, which is $N(j,i)$. Therefore it suffices to show that
$N(i,j)=N(j,i)$.

The cases $i=0$ and $i=1$ correspond to starting at $A_0$ and
$A_1$. {}From Figure \ref{fig-2}, we can check that
 $N(i,j)=N(j,i)$ directly. We have
$$\begin{array}{c| c c c c c}
i\backslash j & 0 & 1 & 2 & 3 & 4   \\ \hline
0      & 1 & 2 & 1 & 0 & 0  \\
1      & 2 & 5 & 9 & 5 & 1  \\
2      & 1 & 9 & * & * & *  \\
3      & 0 & 5 & * & * & *  \\
4      & 0 & 1 & * & * & *
\end{array}$$
and $N(i,j)=0$ if one of $i,j$ is $0$ or $1$ and the other is
great than $4$.

In the case $i\ge 2$, we count the number of $K$-paths from $A_i$
to $B_j$ according to its intersection with the line $y=-2x$. From
Figure \ref{fig-2}, we see that there are $4$ possible
intersection points. We have
\begin{multline}
\qquad N(i,j)=B(2j-i+2,2i-j-4)+3B(2j-i+1,2i-j-2)\\
+3B(2j-i,2i-j)+B(2j-i-1,2i-j+2),\qquad
\end{multline}
where $B(a,b)=\binom{a+b}{b}$.

Let $M(a,b)$ be defined by
$$M(a,b)=B(a+2,b-4)+3B(a+1,b-2)+3B(a,b)+B(a-1,b+2).$$
Then $N(i,j)=M(2j-i,2i-j)$. We need to show that $M(a,b)=M(b,a)$
for $a+b\ge 2$, which implies $N(i,j)=N(j,i)$ for $i\ge 2$. Using
the basic identity of binomial coefficients
$B(c,d)=B(c-1,d)+B(c,d-1)$ for all integers $c$ and $d$, when
$a+b\ge 2$, we have
\begin{multline}
M(a,b)
=B(a-4,b+2)+3B(a-3,b+1)+6B(a-2,b)+7B(a-1,b-1)  \\
+ 6B(a,b-2)+3B(a+1,b-3)+B(a+2,b-4). \label{e-mab}
\end{multline}
In the following Figure \ref{fig-3}, every number we put at a
point is the sum of the numbers at points that are to the left of
it or under it. This corresponds to the formula
$B(c,d)=B(c-1,d)+B(c,d-1)$.

\begin{figure}[ht]
\begin{center}
\input{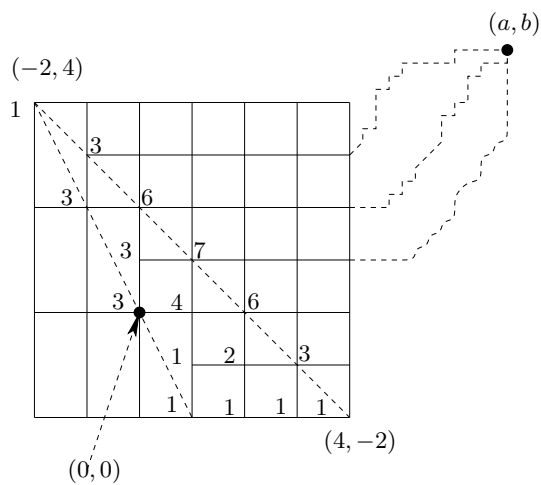}
\end{center}
\caption{Proof of equation \eqref{e-mab} by picture.\label{fig-3}}
\end{figure}

By the symmetry property $B(c,d)=B(d,c)$, we have $M(a,b)=M(b,a)$.
This completes the proof.
\end{proof}

\begin{rem} Observe that the symmetry property of the
numbers $(1,3,6,7,6,3,1)$ along the diagonal in Figure \ref{fig-3}
implies \eqref{e-mab}. A bijective proof of this symmetry will
induce a bijective proof of $N(i,j)=N(j,i)$, and then a bijective
proof of Lemma \ref{l-nij}.\end{rem}

\begin{proof}[Proof of Theorem \ref{t-newg}]
Let $G(m,n)$ be the number of $K$-paths starting at $(-2m,-2m)$ and
ending at
$(2n,2n)$. We will prove that $G(m-1,n+1)=G(m,n)$ for all $m>0$.
Then by induction, $G(m,n)=G(0,m+n)=a_{m+n}$.

We give the bijection as follows. Given a $K$-path $P$ from
$(-2m,-2m)$ to $(2n,2n)$, we separate it by the two lines
$y+2=-2(x+2)$ and $x=-2y$ into three parts: a $V_2$-path $P_1$,
followed by a $K$-path $P_2$, followed by an $H_2$-path $P_3$.

Applying the bijection in the sliding lemma (Lemma \ref{l-nij})
for $P_2$, we get $P^\prime _2$, a $K$-path starting on the line
$y=-2x$, ending on the line $x-2=-2(y-2)$. Then
$P^\prime=P_1P_2^\prime P_3$ with starting point $(-2m+2,-2m+2)$
is the desired $K$-path.

A similar argument gives the inverse bijection.
\end{proof}

This bijective proof of Theorem \ref{t-newg} is not very
desirable, though it is sufficient to prove the determinant
formula \eqref{e-un-res}. The proof relies on the sliding lemma,
whose proof involves a case by case bijection that is not
explicitly given. We would prefer a \emph{natural} bijection for
the sliding lemma that preserves the nonintersecting properties of
$K$-paths. This is because such a bijection would give rise to a
bijection from $\mathcal{U}_K(n)$ to $\mathcal{U}_L(n)$: we could
slide the $n$-tuples of $K$-paths in $\mathcal{U}_K(n)$ so that
all the paths are above the line $x=-2y$. Then the resulting paths
would all be $H_2$-paths that can be easily converted into normal
paths in $\mathcal{U}_L(n)$.

\section{Trinomial Coefficients and $a_n$}

In this section, we introduce $\B$-paths that are counted by
\emph{trinomial coefficients}. The trinomial coefficient $\B(a,b)$
is defined by
$$\B(a,b)= [x^ay^b] (x^2+xy+y^2)^{{a+b\over 2}},$$
if $a+b$ is even, and $\B(a,b)=0$ otherwise.

The trinomial coefficients $T(a,b)$  have a simple combinatorial
interpretation: We call a path $P$ a $\B$-path if each step of $P$
is $(2,0)$ or $(1,1)$ or $(0,2)$. Then the number of $\B$-paths
that start at $(0,0)$ and end at $(a,b)$ is $\B(a,b)$. This
follows easily from the definition of $T(a,b)$. See the following
Figure \ref{fig-4}, in which dots represent vertices of
$\B$-paths.
\begin{figure}[ht]
$$\includegraphics[width=6cm]{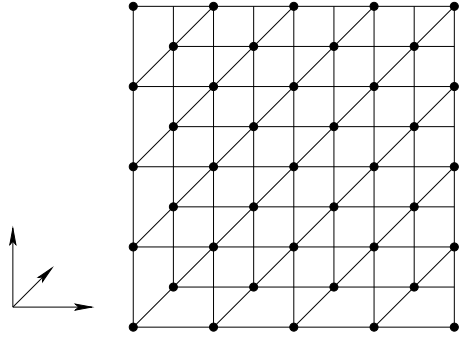}
$$
\caption{The grid for $T$-paths. \label{fig-4}}
\end{figure}

Using $\B$-paths, we can give a better bijective proof of the
sliding lemma. In addition, we find a new determinant identity
\eqref{e-unt}.

Let $n=(a+b)/2$. We can get another formula for $\B(a,b)$ in the
following way. We write $(x^2+xy+y^2)^n=(x(x+y)+y^2)^n$ and use
the binomial theorem twice:
\begin{align*}
(x(x+y)+y^2)^n & = \sum_{k=0}^n {n\choose k}
y^{2(n-k)}x^k\sum_{l=0}^{k}
{k\choose l } x^ly^{k-l} \\
&=\sum_{a=0}^{2n} \left(\sum_{k=0}^a {n\choose k} {k\choose
a-k}\right) x^a y^{2n-a}.
\end{align*}
So
\begin{align}\label{e-t2ab}
T(a,b)= \sum_{k=0}^a {n\choose k} {k\choose a-k}.
\end{align}
This algebraic fact gives another combinatorial explanation of
$\B(a,b)$:
\begin{lem} \label{l-tv}
The number of paths from $(0,-2m)$ to $(i,-i)$, in which the part
below the line $y=-2x$ is a $V_2$ path, and the other part is a
normal path, is equal to the number of $\B$-paths from $(0,-2m)$
to $(i,-i)$, which is $\B(i,2m-i)$.
\end{lem}
\begin{proof}
For a given path $P$ from $(0,-2m)$ to $(i,-i)$, with the part
$P_1$ below the line $y=-2x$ a $V_2$ path, and the other part
$P_2$ a normal path, it is clear that $P_1$ must end at a point
$(j,-2j)$ for some $j\ge 0$, and this $j$ is unique.

We observe that the number of horizontal steps in $P_1$ is $j$,
which equals the total number of  steps in $P_2$. Therefore, we
can associate to each horizontal step in $P_1$ a step in $P_2$,
with order preserved. We call this new path $Q$. Clearly, $Q$ is a
$\B$-path, since each step of $Q$ is a $(0,2)$-step, which is kept
from $P_1$, or a $(1,1)$-step, by associating a vertical step in
$P_2$ to a horizontal step in $P_1$, or a $(2,0)$-step, by
associating a horizontal step in $P_2$ to a horizontal step in
$P_1$. Since the above procedure is a rearrangement of the steps
in $P_1P_2$, $Q$ is a $\B$-path from $(0,-2m)$ to $(i,-i)$. So $Q$
is the desired $\B$-path. The above procedure is clearly
reversible.
\end{proof}

By reflecting in $y=-x$, we get
\begin{lem}\label{l-th}
The number of paths from $(i,-i)$ to $(0,2m)$, in which the part
above the line $x=-2y$ is an $H_2$ path, and the other part is a
normal path, is also $\B(i,2m-i)$.
\end{lem}

Pictures for generalizations of these lemmas can be found in
Figures \ref{fig-5} and \ref{fig-6}. These lemmas correspond to
the case $r=2$.

\begin{thm}\label{t-kte} The number of $K$-paths from $(0,-2m)$ to $(2n,0)$
is equal to the number of $\B$-paths from $(0,-2m)$ to $(2n,0)$,
which is $T(2n,2m)$. Moreover, we have the following identity.
\begin{align}
U_n= \det\left(T(2i,2j)\right)_{0\le i,j \le n-1}.\label{e-unt}
\end{align}
\end{thm}

\begin{proof}[Proof of Theorem \ref{t-kte}]
We can split any $K$-path from $(0,-2m)$ to $(2n,0)$ into two
parts: one  ends at $(i,-i)$ and the other starts at $(i,-i)$ for
some $i\ge 0$ (this $i$ is unique). Then using the two bijections
in Lemmas \ref{l-tv} and \ref{l-th}, we have a bijective proof of
the first part of the corollary.

We have shown in last section that $U_n$ equals the number of
$n$-tuples of nonintersecting $K$-paths from $P_0,\ldots ,P_{n-1}$
to $Q_0,\ldots , Q_{n-1}$, where $P_i=(-2i,-2i)$, and
$Q_i=(2i,2i)$ for $i=0,1,\ldots ,n-1$. It is clear (see Figure
\ref{fig-1}) that it is still true if we replace $P_i$ by
$\tilde{P}_i=(0,-2i)$, and $Q_i$ by $\tilde{Q}_i=(2i,0)$. But from
the first part of this corollary, the number of $K$-paths from
$\tilde{P}_i$ to $\tilde{Q}_j$ is $T(2j,2i)=T(2i,2j)$ for all
$0\le i,j \le$. Then the identity \eqref{e-unt} follows.
\end{proof}
\begin{rem}
The identity \eqref{e-unt} has a generalization in Section
\ref{sec-genKT}. Note that $U_n$ does not equal the number of
$n$-tuples of nonintersecting $\B$-paths from $P'_0,\ldots,
P'_{n-1}$ to $Q'_0,\ldots ,Q'_{n-1}$, because their steps  can
cross without meeting at a vertex of the $\B$-paths.
\end{rem}

\begin{dfn}
We call a path $P$ a $K\B$-path if it satisfies the following
conditions.
\begin{enumerate}
\item[1.] The path $P$ never goes above the diagonal.
\item[2.] The part of $P$ that is to the left of the line $x=0$ is a $V_2$
path.
\item[3.] The part of $P$ in the fourth quadrant is a $\B$-path.
\item[4.] The part of $P$ that is above the line $y=0$ is an $H_2$ path.
\end{enumerate}
\end{dfn}

\begin{thm}\label{t-intgkt}
The number of $K\B$-paths from $(-2m,-2m)$ to $(2n,2n)$ is
$a_{m+n}$ for all $m+n\ge 0$.
\end{thm}
We give three bijective proofs of this theorem. The first
bijective proof establishes the bijection from $K$-paths to
$K\B$-paths. A sliding lemma for $K\B$-paths will then yield a
sliding lemma for $K$-paths. We find that it is much easier to
slide $K\B$-paths: we can slide slowly and we can also slide fast.
We give the fast sliding in our bijection from $K\B$-paths to
$V_2$-paths. This is the second proof. The slow sliding will be
given in Section \ref{sec-genKT} in a more general setting. This
yields the third proof. We suspect that the sliding lemma for
$K$-paths resulting from our second and third bijections are
natural, i.e., preserve the nonintersecting property.

\begin{proof}[Bijection from $K\B$-paths to $K$-paths]
We first uniquely separate $P$, according to its intersections
with the lines $x=0$ and $y=0$,  into three parts, a $V_2$-path
$P_1$, followed by a $\B$-path $P_2$, followed by an $H_2$-path
$P_3$, such that $P_1$ ends with a horizontal step and $P_3$
starts with a vertical step, except that $P_1$ and $P_3$ may be
empty.

{}From Theorem \ref{t-kte}, we can get a $K$-path $P_2^\prime$ from
$P_2$ without changing the starting and ending points. Then
$P^\prime=P_1P_2^\prime P_3$ is the desired $K$-path. This
procedure is clearly reversible.
\end{proof}

The next proof relies highly on Lemmas \ref{l-tv} and \ref{l-th}.
The bijection $\phi_v$ for Lemma \ref{l-tv} maps a $\B$-path $P$
to a $V_2$-path $P_V$ followed by a normal path $P_N$, in which
the number of horizontal steps in $P_V$ equals the total number of
steps in $P_N$. Given the starting point $S(P)$ and ending point
$E(P)$ of $P$, we can predict the position of $E(P_V)=S(P_N)$:
$E(P_V)$ must lie on the line with slope $-2$ and passing through
the point $O$, which is determined by the conditions that
$O\rightarrow S(P)$ is vertical and the slope of $O \rightarrow
E(P)$ is $-1$.

Similarly the bijection $\phi_h$ for Lemma \ref{l-th} maps a
$\B$-path $P$ to a normal path $P_N$ followed by an $H_2$-path
$P_H$ with similar properties.

\begin{proof}[Fast sliding bijection from $K\B$-paths to $V_2$-paths]
Let $P$ be a $K\B$-path from $(-2m,-2m) $ to $(2n,2n)$ with
$m,n\ge 0$.

We first uniquely separate $P$, according to its intersections
with the lines $x=0$ and $y=0$,  into a $V_2$-path $P_1$, followed
by a $\B$-path $P_2$, followed by an $H_2$-path $P_3$, such that
$P_1$ ends with a horizontal step and $P_3$ starts with a vertical
step, except that $P_1$ and $P_3$ may be empty. We will map
$P_2P_3$ to a $V_2$ path with the same starting and ending points.

Suppose $S(P_2)=(0,-i)$ and $E(P_2)=(j,0)$. Obviously we can
assume $j>0$ for otherwise $P_2$ is a $V_2$-path and $P_3$ is the
empty path. Applying $\phi_v$ to $P_2$ gives us a $V_2$-path
$P_{2V}$ followed by a normal path $P_{2N}$ with
$S(P_{2N})=(a,-2a+j)$ for some $a>0$. See Figure \ref{fig-n2},
where we did not draw the paths explicitly.

\begin{figure}[ht]
\begin{center}
\input{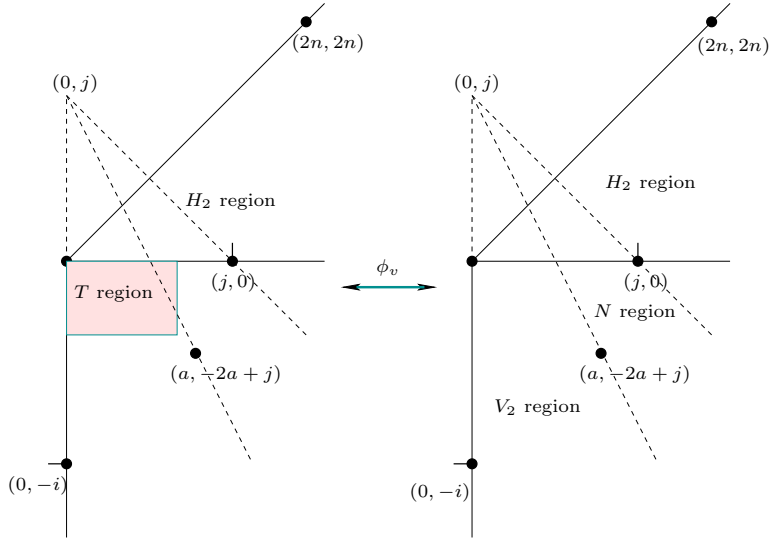}
\end{center}
 \caption{First step of the fast sliding. \label{fig-n2}}
\end{figure}

Draw a vertical line at $S(P_{2N})$, which intersects the diagonal
at $(a,a)$. It is easy to check that the the total number of steps
of $P_{2N}$ is $a$. Factor the $H_2$-path $P_3$, according to its
intersection with the line $y=a$, into $P_{31}P'_{3}$ such that
$P'_{3}$ starts with a vertical step. The the number of vertical
steps of $P_{31}$ equals $a$. Applying $\phi_h^{-1}$ to
$P_{2N}P_{31}$ gives us a $\B$-path $P'_2$. See Figure
\ref{fig-n3}.

\begin{figure}[ht]
\begin{center}
\input{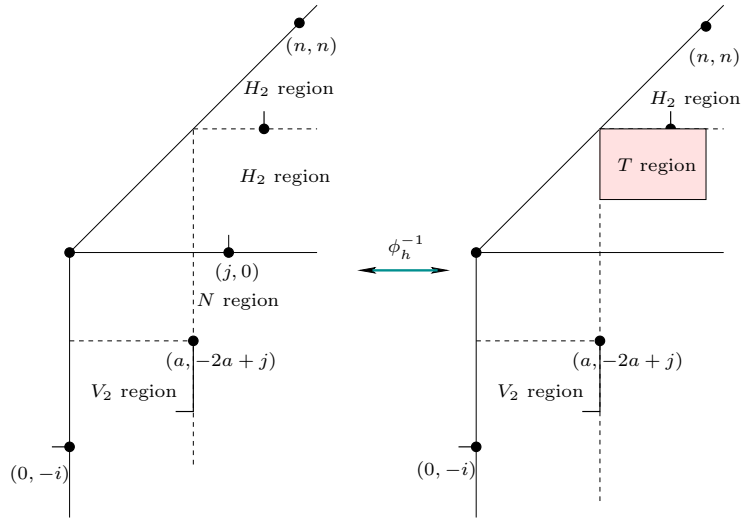}
\end{center}
 \caption{Second step of the fast sliding.\label{fig-n3}}
\end{figure}

Now we are left to map $P_{2V}P'_2 P'_3$ to a $V_2$-path. If we
slide down the path by $(a,a)$, then we met the same situation as
for the path $P_2P_3$. Repeat the above procedure we can finally
obtain the desired $V_2$-path. The procedure is reversible because
of the required
 conditions of
ending with a horizontal step or starting with a vertical step, as
shown in Figures \ref{fig-n2} and \ref{fig-n3}.
\end{proof}

\section{Generalizations of $K$-paths and $KT$-paths \label{sec-genKT}}
Let $\rr{g}_n=\frac{1}{rn+1}{(r+1)n\choose n}$ be the number of
$r+1$-ary trees with $n$ nodes, and
$$\rr{g}(x)=\sum_{n\ge 0} \rr{g}_n x^n$$
be the generating function. Then $\rr{g}(x)$ satisfies the
following functional equation.
$$\rr{g}(x)=1+x\left(\rr{g}(x) \right)^{r+1}.$$

For $r=1$, $g_n^{(1)}$ is the Catalan number. It is well-known
that the Hankel determinants of the Catalan generating function
are all $1$. We have studied the the case $r=2$. We wish to say
something about the Hankel determinants of $\rr{g}(x)$ for $r\ge
3$.

Since $H_n(\rr{g}(x))$ does not factor for $r\ge 3$, a formula
like \eqref{e-typea} is unlikely. However, we find generalizations
of \eqref{e-un-res} (which is the same as \eqref{e-un}),
\eqref{e-matrix-identity}, and \eqref{e-unt}. They are given by
\eqref{e-grij}, \eqref{e-98}, and \eqref{e-99}. Their algebraic
proofs can be found in Section \ref{sec-alg-pf}. Except for
$H_n(\rr{g}(x))$, we do not have nonintersecting paths
interpretation of these determinants.

We have natural generalizations of the concepts in the last two
sections. A path is an $H_r$ path if each step is either $(r,0)$
or $(0,1)$. Similarly, a path is a $V_r$ path if each step is
either $(1,0)$ or $(0,r)$. The following is equivalent to a
special case of a classical result given (without proof) by
Barbier \cite{barbier}. (See \cite{gs} for a new proof and further
references.)

\begin{prop}\label{p-grn}
The number of $V_r$-paths $($or $H_r$-paths$)$ from $(0,0)$ to
$(rn,rn)$ that never go above the  diagonal is $\rr{g}_n$.
\end{prop}

A path is a $\rr{\B}$-path if each of its step is $(r,0)$,
$(r-1,1)$,
 \ldots , or $(0,r)$.
For any path $P$, we denote by $S(P)$ the starting point, $E(P)$
the ending point, and $L(P)$ the number of steps in $P$.

For a $\rr{\B}$-path $P$ with $S(P)=(0,0)$ and $L(P)=k$, $E(P)$
must lie on the line $y=-x+rk$. So to compute $L(P)$, we take the
sum of the $x$-coordinate and $y$-coordinate of $E(P)-S(P)$, and
divide by $r$. A normal path is also a $\B^{(1)}$-path, a $V_1$
path, and an $H_1$ path, and a $\B$-path is a $\B^{(2)}$-path.

Let $\rr{\B}(a,b)$ be the number of $\rr{\B}$-paths from $(0,0)$
to $(a,b)$. Then $\rr{\B}(a,b)=0$ if $a+b$ is not divisible by
$r$, so we can suppose $a=ri-s$ and $b=rj+s$ for some $i,j$ and
$0\le s \le r-1$. We have
$$\rr{\B}(ri+s,rj-s) = [x^{ri+s}y^{rj-s}] (x^r+x^{r-1}y+\cdots
+y^r)^{i+j}.$$ Since the right-hand side of the above equation is
homogeneous in $x$ and $y$, we can write it in terms of one
variable $t$, where $t=y/x$.

Let $\alpha=1+t+t^2+\cdots +t^r$ and
$\beta=\alpha/t^{r}=1+t^{-1}+\cdots t^{-r}$. Then
\begin{align}
\label{e-CT-trinomial} T(ri+s,rj-s)=\ct t^s \alpha^i\beta^j=\ct
t^{s-rj} \alpha^{i+j},
\end{align}
where $\ct$ means to take the constant term of a Laurent
polynomial of $t$.

\def\tt{\mathbf{ \B}}
\def\vv{\mathbf{ V}}
\def\hh{\mathbf{ H}}

\begin{dfn}
We call a path $P$ a $\rr{K}$-path if it satisfies the following
conditions.
\begin{enumerate}
\item[1.] The path $P$ never goes above the diagonal.
\item[2.] The part of $P$ that is to the left of the line $y=-rx$ is a $V_r$
path.
\item[3.] The part of $P$ between the two lines $y=-rx$ and $x=-ry$ is a
$\B^{(r-1)}$-path.
\item[4.] The part of $P$ that is above the line $x=-ry$ is an $H_r$ path.
\end{enumerate}
\end{dfn}

\begin{dfn}
We call a path $P$ a $\rr{K\B}$-path if it satisfies the following
conditions.
\begin{enumerate}
\item[1.] The path $P$ never goes above the diagonal.
\item[2.] The part of $P$ that is to the left of the line $x=0$ is a $V_r$
path.
\item[3.] The part of $P$ in the fourth quadrant is a $\B^{(r)}$-path.
\item[4.] The part of $P$ that is above the line $y=0$ is an $H_r$ path.
\end{enumerate}
\end{dfn}
For example, a $K$-path is a $K^{(2)}$-path, and a $K\B$-path is a
$K\B^{(2)}$-path.

Let $\mathbb{K}(m,n,r)$ be the set of all $\rr{K}$-paths from
$(-mr,-mr)$ to $(nr,nr)$. Let $\mathbb{T}(m,n,r,s)$ be the set of
all $\rr{K\B}$-paths from $(s-mr,s-mr)$ to $(nr+s,nr+s)$. Now we
can state our main results.
\begin{thm}\label{t-new-grnk}
The cardinality of $\mathbb{K}(m,n,r)$ is $\rr{g}_{n+m}$ for all
$m$ and $n$. The cardinality of $\mathbb{T}(m,n,r,s)$ is also
$\rr{g}_{n+m}$ for all $m,n$ and $s$.
\end{thm}
As in the case $r=2$, if $m\le 0$ $($or $n\le 0$$)$, then
$\rr{K}$-paths and $\rr{K\B}$-paths are in fact $H_r$-paths (or
$V_r$-paths), and in these cases, Theorem \ref{t-new-grnk} follows
from Proposition \ref{p-grn}. The idea of the proof of Theorem
\ref{t-new-grnk} is to show that
$|\mathbb{T}(m,n,r,s)|=|\mathbb{T}(m,n,r,s-1)|$ for all $1\le s\le
r$. Then $|\mathbb{T}(m,n,r,s)|=\rr{g}_{m+n}$ follows by
induction. We will give a bijection from $\mathbb{T}(m,n,r,0)$ to
$\mathbb{K}(m,n,r)$.

The bijective proof we are going to give highly relies on the
following lemma, especially on the bijection from $\tt(i,j)$ to
$\vv(i,j)$, which is a generalization of Lemma \ref{l-tv}.

\begin{lem}\label{l-gtv}
The following four sets  all have cardinality $\B^{(r)}(j,ri-j)$.
\begin{enumerate}
\item The set $\tt (i,j)$ of all $T^{(r)}$-paths from $(0,-ri)$ to $(j,-j)$.
\item The set $\tt '(i,j)$ of all $T^{(r)}$-paths from $(j,-j)$ to
$(ri,0)$.
\item The set $\vv(i,j)$ of all paths from $(0,-ri)$ to $(j,-j)$, with the part
before the line $y=-rx$ a $V_r$ path, and the part after the line
$y=-rx$ a $\B^{(r-1)}$-path.
\item The set $\hh(i,j)$ of all  paths from $(j,-j)$ to $(ri,0)$ with the part
before  the line $x=-ry$ a $\B^{(r-1)}$-path, and the part after
the line $x=-ry$ an $H_r$ path.
\end{enumerate}
\end{lem}

\begin{proof}
We construct only the bijection from $\tt(i,j)$ to $\vv(i,j)$. The
bijection from $\tt '(i,j)$ to $\hh(i,j)$ is similar. The
bijection from $\vv(i,j)$ to $\hh(i,j)$ and the bijection from
$\tt(i,j)$ to $\tt'(i,j)$ are given by reflecting in the line
$y=-x$.

For any given path $T\in \tt(i,j)$ the steps in $T$ are $(r-k,k)$
for $k=0,1,\ldots ,r$. We first replace all the steps in $T$ that
are not $(0,r)$ with steps $(r,0)$. Then we get a path $T_1$, with
$E(T_1)$ on the line $y=-x$. Changing every $(r,0)$-step in $T_1$
into $(1,0)$ will give us a $V_r$ path $V_1$ with $E(V_1)$
 on the line $y=-rx$. By removing all the $(0,r)$-steps in $T$, and changing
every $(r-k,k)$-step into $(r-k-1,k)$ for $k=0,1,\ldots r-1$, we
get a $\B^{(r-1)}$-path $V_2$. Then $V=V_1V_2\in \vv(i,j)$ is the
desired path. One way to see that $E(V)=(j,-j)$ is that if we
regard an $(r-k,k)$-step as a $(1,0)$ step followed by a
$(r-k-1,k)$-step for $k=0,1,\ldots ,r-1$, then the bijection we
gave is just a rearrangement of the steps in $T$.

The inverse procedure is as follows. For a given $V'\in \vv(i,j)$,
the line $y=-rx$ divides $V'$ into a $V_r$ path $V'_1$ followed by
a $\B^{(r-1)}$ path $V'_2$. Suppose $E(V'_1)=(k',-rk')$ for some
$k'$. We can see that the number of $(1,0)$-steps in $V'_1$, which
is $k'$, equals the total number of steps in $V'_2$, which is
${\left((j-k')+(-j+rk')\right)}/{(r-1)}=k'$. Then we can associate
to each $(1,0)$ step in $V_1'$ a step in $V_2'$, with order
preserved. This gives us a $\rr{\B}$-path $T'\in \tt(i,j)$. In
Figure \ref{fig-5}, we give an example of the case $r=4$.
\begin{figure}[ht]
\begin{center}
\input{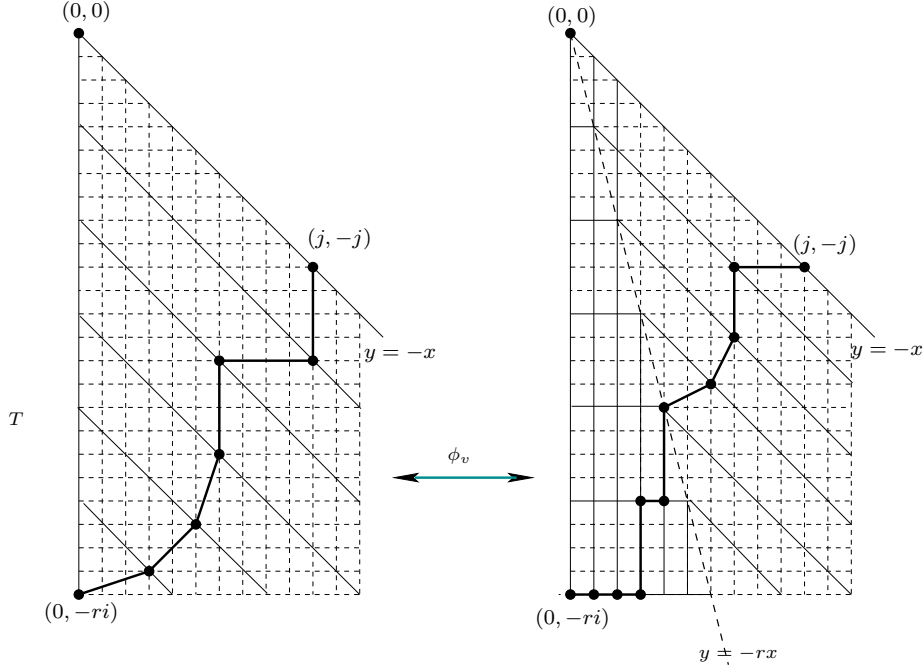}
\end{center}
\caption{A $\rr{\B}$-path $T$ and its image under
$\phi_v$.\label{fig-5}}
\end{figure}
The above two procedures are clearly inverse to each other.
\end{proof}

\begin{rem}
The bijection from $\tt(i,j)$ to $\vv(i,j)$ we gave originated
from the following algebraic fact.

For $r\ge 3$, there are many ways to group and expand the
polynomial $(x^r+x^{r-1}y+\cdots +y^r)^{m}$. We find the following
way has a nice combinatorial explanation.

\begin{align*}
(x(x^{r-1}+x^{r-2}y+\cdots y^{r-1})+y^r)^m &=
\sum_{j=0}^m {m\choose j} y^{r(m-j)} x^j(x^{r-1}+x^{r-2}y+\cdots y^{r-1})^j \\
&=\sum_{i=0}^{rm} \left(\sum_{j=0}^{m} {m\choose j}
\B^{(r-1)}(rj-i, i-j)\right)x^{i} y^{rm-i}
\end{align*}
So we have
\begin{align}\label{e-trexpand}
T^{(r)}(i,rm-i)=\sum_{j=0}^{m} {m\choose j} \B^{(r-1)}(rj-i, i-j).
\end{align}

In Figure \ref{fig-5}, $|\vv(i,j)|$ can be counted according to
the intersection points of the paths in $\vv(i,j)$ with the line
$y=-rx$. This yields \eqref{e-trexpand}.
\end{rem}

\begin{figure}[ht]
\begin{center}
\input{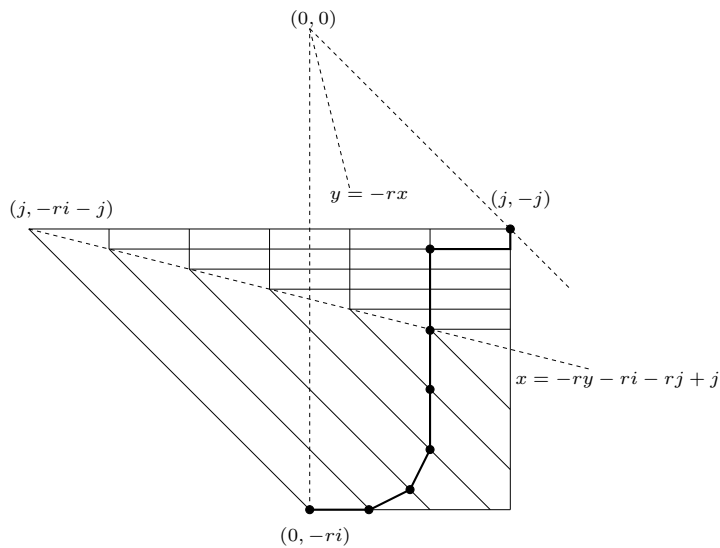}
\end{center}
\caption{The image of $T$ under $\phi_h$. \label{fig-6}}
\end{figure}

We denote the bijection from $\tt(i,j)$ to $\vv(i,j)$ by $\phi_v$,
and the  bijection from $\tt '(i,j)$ to $\hh(i,j)$ by $\phi_h$.
One thing we should mention is that neither $\phi_h$ nor $\phi_v$
changes the starting point and the ending point. The path in
Figure \ref{fig-6} is obtained from the $\rr{\B}$-path $T$ in
Figure \ref{fig-5} by applying $\phi_h$.

Applying $\phi_v$ to a $\rr{\B}$ path will give us a $V_r$ path
followed by a $\B^{(r-1)}$ path, in which the number of horizontal
steps in the $V_r$ path equals the total number of  steps in the
$\B^{(r-1)}$ path. We can locate the ending point of the $V_r$
path by the following three easy steps. (See Figure \ref{fig-5}.)
\begin{enumerate}
\item Draw a vertical line at $S(T)$.
\item Draw a line of slope $-1$ at $E(T)$.
\item At the intersecting point of the above two lines, draw a line of slope
$-r$. Then this is the line on which the ending point of the $V_r$
path must lie.
\end{enumerate}
\def\blv{\mbox{BL}_v}
\def\blh{\mbox{BL}_h}
We call the line obtained in the above three steps the
\emph{bisecting line} $\blv(T)$ of $\phi_v(T)$. For any $\rr{T}$
path $T$, with $S(T)=(0,-ri)$ and $E(T)=(j,-j)$,
 $\blv(T)$ is
$y=-rx$, which is independent of $i$ and $j$.

A similar argument for $\phi_h$ holds for a $\rr{T}$-path $T'$.
The corresponding three steps are given as follows. (See Figure
\ref{fig-6}.)
\begin{enumerate}
\item Draw a horizontal line at $E(T')$.
\item Draw a line of slope $-1$ at $S(T')$.
\item At the intersecting point of the above two lines, draw a line of slope
 $-1/r$.
Then this is the line on which the starting point of the
$H_r$-path must lie.
\end{enumerate}
We call the line obtained in the above three steps the bisecting
line $\blh(T')$ of $\phi_h(T')$.

The following is a generalization of Theorem \ref{t-kte}.
\begin{lem}\label{l-94}
The number of each of the following four kinds of paths from
$(0,-ri+s)$ to $(rj+s,0)$   equals to $T(rj+s,ri-s)$.
\begin{enumerate}
\item All $\rr{\B}$-paths from $(0,-ri+s)$ to
$(rj+s,0)$.
\item All paths from $(0,-ri+s)$ to
$(rj+s,0)$ consisting of a $V_r$ path, followed by a $T^{(r-1)}$
path, in which the number of $(1,0)$-steps in the $V_r$ path
equals the total number of steps in the $T^{(r-1)}$ path.
\item All the paths from $(0,-ri+s)$ to
$(rj+s,0)$  consisting of a $T^{(r-1)}$ path, followed by an $H_r$
path, in which the number of $(0,1)$-steps in the $H_r$ path
equals the total number of steps in the $T^{(r-1)}$ path.
\item All the paths from $(0,-ri+s)$ to
$(rj+s,0)$, with the part before the line $y=-rx+s$ a $V_r$ path,
between the lines $y=-rx+s$ and $x=-ry+s$ a $T^{(r-1)}$ path, and
the part after the line $x=-ry+s$ an $H_r$ path.
\end{enumerate}
\end{lem}
\begin{proof}
Part $1$ follows from the definition of $\B^{(r)}(rj+s,ri-s)$.
Part $2$ and part $3$ are obvious by Lemma \ref{l-gtv}, so we need
only prove part $4$.

For a given $\B^{(r)}$ path $T$ from $(0,-ri+s)$ to $(rj+s,0)$, we
can uniquely separate it by the line $y=-x+s$ into a $\B^{(r)}$
path $T_1$ followed by a $\B^{(r)}$ path $T_2$. Applying $\phi_v$
on $T_1$, we get a $V_r$ path $T_{1V}$ followed by a $\B^{(r-1)}$
path $T_{1T}$. Applying $\phi_h$ on $T_2$, we get a $\B^{(r-1)}$
path $T_{2T}$ followed by an $H_r$ path $T_{2H}$. Using the three
steps for finding $\blv(T_1)$, we see that  $E(T_{1V})$ must lie
on the line $y=-rx+s$, since the line $y=-x+s$ intersects the line
$x=0$ at $(0,s)$. Similarly, $S(T_{2T})$ must lie on the line
$x=-ry+s$, since the line $y=-x+s$ intersect the line $y=0$ at
$(s,0)$. Together with the fact that
$E(T_{1T})=E(T_1)=S(T_2)=S(T_{2T})$, we see that $T_{1T}T_{2T}$ is
also a $T^{(r-1)}$ path and the path $T_{1V}T_{1T}T_{2T}T_{2H}$ is
the desired path.

The above procedure is clearly invertible.
\end{proof}

The bijection in the above proof will induce a bijection from
$\mathbb{T}(m,n,r,0)$ to $\mathbb{K}(m,n,r)$. We will see this in
the proof of Theorem \ref{t-new-grnk}.

For any $1\le s \le r$, let $\pa(i,j)$ be the set of $\bk$-paths
from $(0,-(ri-s))$ to $(rj+s,1)$, and $\pb(i,j)$ be the set of
$\bk$-paths from $(-1,-(ri-s+1)$ to $(rj+s-1,0)$. Then we have the
following lemma, which will induce the bijection from
$\mathbb{T}(m,n,r,s)$ to $\mathbb{T}(m,n,r,s-1)$.

\begin{lem}[Slow Sliding Lemma]\label{l-mainlemma}
There is a bijection from $\pa(i,j)$ to $\pb(i,j)$ for all $i,j$.
\end{lem}
We will give two proofs for this lemma. The algebraic proof will
be given in the next section. The bijective proof is as follows.

\begin{proof}[Bijective proof of Lemma \ref{l-mainlemma}]
For any given $P\in \pa (i,j)$, we uniquely factor $P$, according
to its intersections with the lines $y=0$ and $y=-x+s$, into
$P_1P_2P_3$, where we require $P_3$ to start with a vertical step.
In the left picture of Figure \ref{fig-7a}, we marked each
intersection point by a $\circ$.

\begin{figure}[ht]
\begin{center}
\input{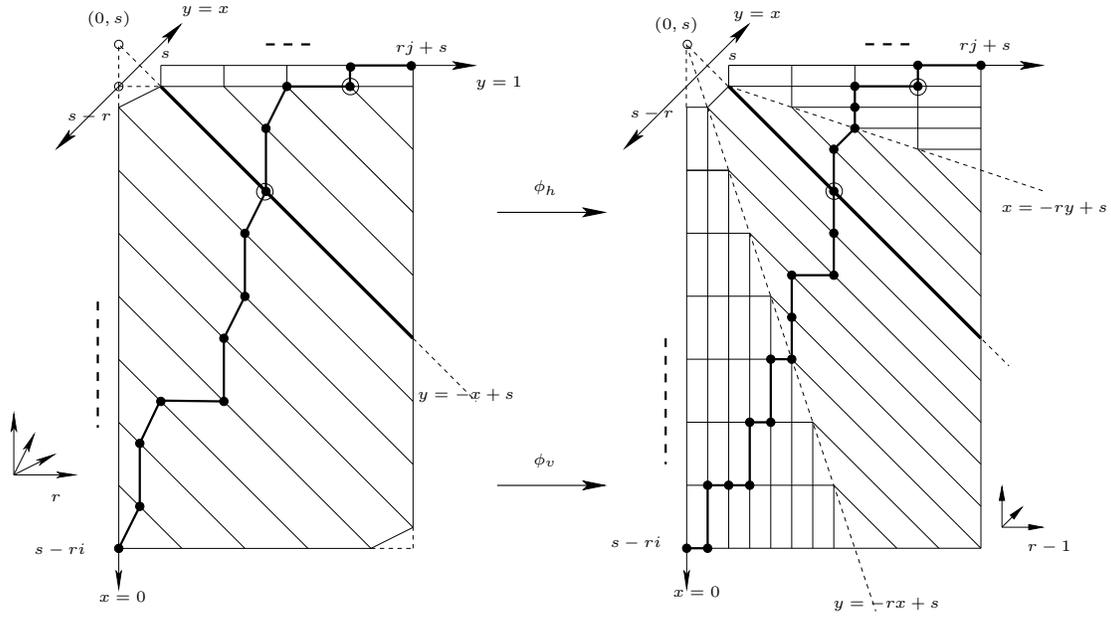}
\end{center}
\caption{First step of the slow sliding bijection.\label{fig-7a}}
\end{figure}

Now we apply $\phi_v$ to $P_1$ to obtain a $V_r$-path followed by
a $T^{(r-1)}$-path, and apply $\phi_h$ to $P_2$ to obtain a
$T^{(r-1)}$-path followed by an $H_r$-path. The bisection lines
are $y=-rx+s$ and $x=-ry+s$, as drawn in the right picture of
Figure \ref{fig-7a}.

Denote by $P'$ the whole path obtained this way. We uniquely
factor $P'$, according to its intersections with the lines $x=1$
and $y=-x+s-r+1$ into $P_1'P_2'P_3'$, where we require $P_1'$ to
end with a horizontal step. In the left picture of Figure
\ref{fig-7b}, we marked each intersection point by a $\Box$.

\begin{figure}[ht]
\begin{center}
\input{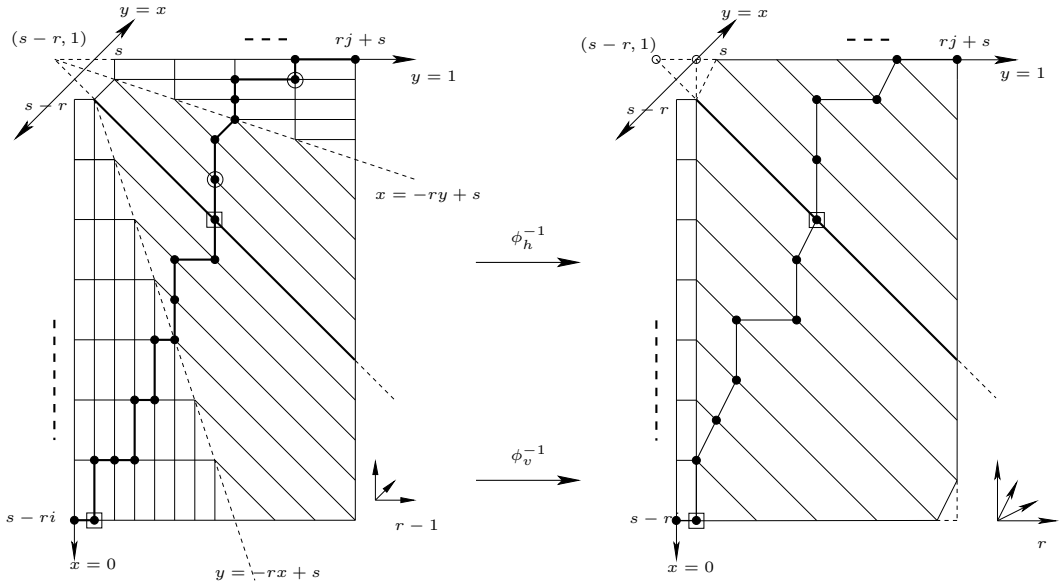}
\end{center}
\caption{Second step of the slow sliding bijection.\label{fig-7b}}
\end{figure}

Now we apply $\phi_v^{-1}$ to $P_2'$ to obtain a $T^{(r)}$-path
$Q_2$, and apply $\phi_h^{-1}$ to $P'_3$ to obtain a
$T^{(r)}$-path $Q_3$. See the right picture of Figure
\ref{fig-7b}. We need to check that the three lines $y=-x+s-r+1$,
$x=-ry+s$, and $y=1$ intersect at the point $(s-r,1)$, so that
$\phi_h^{-1} (P'_3)$ is well defined.

Finally, let $Q$ be obtained from $P_1' Q_2Q_3$ by sliding down by
$(1,1)$. Then $Q$ is the desired path. Every step in the above
procedure is invertible. This completes the proof.
\end{proof}

\begin{proof}[Proof of Theorem \ref{t-new-grnk}]
First we construct the bijection from $\mathbb{T}(m,n,r,0)$ to
$\mathbb{K}(m,n,r)$. This shows that
$|\mathbb{T}(m,n,r,0)|=|\mathbb{K}(m,n,r)|$.

  Recall that any $P\in
\mathbb{T}(m,n,r,0)$ is a \ktrpath from $(-rm,-rm)$ to $(rn,rn)$,
and any $Q\in \mathbb{K}(m,n,r)$ is a \krpath with the same ending
points. We can uniquely factor $P$, according to its intersections
with the lines $x=0$ and $y=0$, into $P_1P_2P_3$, such that $P_1$
is a $V_r$ path ending with a horizonal step and $P_3$ is an $H_r$
path starting with a vertical step, except that $P_1$ and $P_3$
may be empty. Applying the bijection of Lemma \ref{l-94} part $4$
to $P_2$, we get a \krpath $Q_2 =Q_{2a}Q_{2b}Q_{2c}$ from $S(P_2)$
to $E(P_2)$, with $Q_{2a}$ a $V_r$ path ending on the line
$y=-rx$, $Q_{2b}$ a $T^{(r-1)}$ path ending on the line $x=-ry$,
and $Q_{2c}$ an $H_r$ path. Then $P_1Q_{2a}$ is  a $V_r$ path, and
$Q_{2c}P_3$ is an $H_r$ path. So
$Q=P_1Q_2P_3=(P_1Q_{2a})Q_{2b}(Q_{2c}P_3)$ belonging to
$\mathbb{K}(m,n,r)$ is the desired path. The above procedure is
clearly reversible.

Next we construct the bijection from $\mathbb{T}(m,n,r,s)$ to
$\mathbb{T}(m,n,r,s-1)$ for $1\le s \le r$, which implies
$|\mathbb{T}(m,n,r,s)|=|\mathbb{T}(m,n,r,s-1)|$. Any $P\in
\mathbb{T}(m,n,r,s)$ can be uniquely factored, according to its
intersections with the lines $x=0$ and $y=1$, into $P_1P_2P_3$,
such that $P_1$ is a $V_r$ path ending with a horizonal step and
$P_3$ is an $H_r$ path starting with a vertical step, except that
$P_1$ and $P_3$ may be empty. Then $P_2$ is a $KT^{(r)}$-path.

Applying the bijection of Lemma \ref{l-mainlemma} to $P_2$, we get
a \ktrpath $Q_2=Q_{2a}Q_{2b}$, with $Q_{2a}$ a $V_r$ path starting
on the line $x=-1$ and ending on the line $x=0$, and $Q_{2b}$ a
$T_r$ path ending on the line $y=0$. Let $Q_1$ be obtained from
$P_1$ by sliding down by $(1,1)$, and $Q_3$ be obtained from $P_3$
by sliding down by $(1,1)$. Then
$Q=Q_1Q_2Q_3=(Q_1Q_{2a})Q_{2b}Q_3$ is a path from
$(s-1-mr,s-1-mr)$ to $(s-1+nr,s-1,nr)$ that never goes above the
diagonal and with the part before the line $x=0$ a $V_r$ path, the
part between the two lines $x=0$ and $y=0$ a $T_r$ path, and the
part after the line $y=0$ an $H_r$ path. Hence $Q\in
\mathbb{T}(m,n,r,s-1)$ is the desired path. The above procedure is
clearly reversible.

Finally, we use induction to conclude the theorem. By the second
part, it is easy to see that
$|\mathbb{T}(m,n,r,s)|=|\mathbb{T}(m,n,r,r)|$ for $0\le s\le r$.
But $\mathbb{T}(m,n,r,r)$ is in fact $\mathbb{T}(m-1,n+1,r,0)$.
Theorem \ref{t-new-grnk} then follows by induction and the fact
that $|\mathbb{T}(0,m+n,r,0)|=\rr{g}_{m+n}$.
\end{proof}
\begin{rem}
We can also give a fast sliding bijection from
$\mathbb{T}(m,n,r,0)$ to $\mathbb{T}(0,m+n,r,0)$.
\end{rem}

Recall that $\rr{g}(x)$ is the generating function of $r+1$-ary
trees. Let $\rr{f}(x)=\rr{g}(x)-1$. Then $\rr{f}$ satisfies the
following functional equation.
$$\rr{f}=x(1+\rr{f})^{r+1}.$$
If we count $K^{(r)}$-paths from $(-mr,-mr)$ to $(nr,nr)$
according to their intersections with the lines $y=-rx$ and
$x=-ry$, we see that Theorem \ref{t-new-grnk} yields the matrix
identity \begin{multline}\label{e-grij}
\left(\rr{g}_{i+j}\right)\ijn =\left([x^i] \rr{g} (\rr{f})^j
\right)\ijn \\ \left(T^{(r-1)}(rj-i,ri-j)\right)\ijn \left( [x^j]
\rr{g} (\rr{f})^i\right)\ijn,
\end{multline}
where
$$[x^i] \rr{g} (\rr{f})^j=\frac{(r+1)j+1}{(r+1)i+1}{(r+1)i+1\choose
i-j}$$ is the number of $V_r$-paths from $(-ri,-ri)$ to $(j,-rj)$
that never go above the diagonal.

Since the transformation matrices in \eqref{e-grij} are upper (or
lower) triangular matrices with diagonal entries $1$, we have

\begin{align}\label{e-98}
\det\left(\rr{g}_{i+j}\right)\ijn =\det
\left((T^{(r-1)}(rj-i,ri-j)\right)\ijn
\end{align}

A similar argument gives

\begin{align}\label{e-99}
\det\left(\rr{g}_{i+j}\right)\ijn =\det
\left(\rr{T}(ri,rj)\right)\ijn =\det
\left(\rr{T}(ri-s,rj+s)\right)\ijn,\end{align}
 for any $0\le s \le r-1$.

\section{The Algebraic Proof\label{sec-alg-pf}}

The scheme of our algebraic proof of Lemma \ref{l-mainlemma} is by
first representing our object as the constant term of a rational
function, and then evaluating the constant term. This technique is
well-known. For instance, Egorychev \cite{ego} gave many
applications for evaluating combinatorial sums in the context of
residues (equivalent to constant terms). We also use this method,
together the method in Section 4, to give algebraic proofs of
\eqref{e-98} and \eqref{e-99}.

The only thing we need here is the following Proposition
\ref{p-ira}. Its proof is included since the idea of the proof
applies to most of our examples. We will give a different
algebraic proof of
 equations \eqref{e-gf1} and \eqref{e-r}.

Let $B(x,y,t)\in \mathbb{Q}[t,t^{-1}][[x,y]]$. Then $B(x,y,t)$ can
be written as
$$B(x,y,t) = \sum_{i,j=0}^\infty b_{ij}(t)x^iy^j,$$
where $b_{ij}(t)\in \mathbb{Q}[t,t^{-1}]$. Define
$$\ct B(x,y,t) =\sum_{i,j=0}^\infty (\ct b_{ij}(t)) x^iy^j ,$$
where $\ct b_{ij}(t)$ is the constant term of the Laurent
polynomial $b_{ij}(t)$ in $t$.

\def\CT{\mathop{\rm CT}}
\def\tt{t^{-1}}
\def\QQ{\mathbb{Q}}

The general problem in this section is to find the constant term
of the function $(1-P(\tt)x)^{-1}(1-Q(t)y)^{-1}$, for some
specific $P(t), Q(t)\in \mathbb{Q}[t]$.

\begin{prop}\label{p-ira}
Let $P(t)$ and $Q(t)$ be polynomials in $t$, and let $a_{mn}=\CT
P(t)^m Q(t^{-1})^n$. Then
$$\sum_{m,n}a_{mn}x^m y^n$$ is a rational function in $x$ and $y$.
\end{prop}
\begin{proof}
We show that
$$\CT {1\over (1-P(t)x)(1-Q(\tt)y)}$$ is rational, where we work in
the ring $\QQ[t,\tt][[x,y]]$.
 We may assume that $P(t)$ has degree at at least 1. Let $d$ be the
degree of $Q(t)$. Let
$$F= {1\over (1-P(t)x)(1-Q(\tt)y)}={t^d\over (1-P(t)x)(t^d-t^d
Q(\tt)y)}.$$

Since $t^d Q(\tt)$ is a polynomial in $t$ of degree at most $d$,
and the degree of $P(t)$ is at least 1, $F$ has a partial fraction
expansion in $t$ that may be written
\begin{align}
F={1\over R(x,y)} \left( {A(x,y,t)\over 1-P(t)x} +{B(x,y,t)\over
t^d-t^dQ(\tt)y} \right) = {1\over R(x,y)}\left( {A(x,y,t)\over
1-P(t)x} +{B(x,y,t)t^{-d}\over 1-Q(\tt)y} \right) \label{e-ira1}
\end{align}
where $R(x,y)$ is a polynomial in $x$ and $y$, $A(x,y,t)$ and
$B(x,y,t)$ are polynomials in $x$,  $y$, and $t$, and the degree
of $B$ in $t$ is less than $d$.

Now the constant term in $t$ of $B(x,y,t)t^{-d}/( 1-Q(\tt)y)$ is 0
and the constant term in $A(x,y,t)/( 1-P(t)x)$ is
$A(x,y,0)/(1-P(0)x)$. We would like to conclude that
\begin{align}
\CT F={A(x,y,0)\over R(x,y)(1-P(0)x)}.\label{e-ira2}
\end{align}
However, we don't know that $1/R(x,y)$ has a power series
expansion. To avoid this problem, we multiply \eqref{e-ira1} by
$R(x,y)$ to get
$$R(x,y)F=
{A(x,y,t)\over 1-P(t)x} +{B(x,y,t)t^{-d}\over 1-Q(\tt)y}.
$$

Then
$$\CT R(x,y)F={A(x,y,0)\over 1-P(0) x},$$
but since $\CT R(x,y) F= R(x,y)\CT F$, \eqref{e-ira2} follows.
\end{proof}

Since the main idea of this proof is a partial fraction
decomposition, we call this method the \emph{partial fraction
method}. In the following examples, we use formula
\eqref{e-CT-trinomial}. Let $\alpha=1+t+t^2+\cdots +t^r$ and
$\beta=\alpha/t^{r}=1+t^{-1}+\cdots t^{-r}$. Then
$$T(ri+s,rj-s)=\ct t^s \alpha^i\beta^j. $$
In particular,
$${m+n\choose m}=[t^n] (1+t)^{m+n} = \ct (1+t)^m (1+\frac{1}{t})^n.$$

\begin{exa}\label{ex-1}
A different proof of identity \eqref{e-gf1}
\begin{equation*}
\frac{1 - xy}{1 - xy^2 - 3xy - x^2y}= \sum_{i,j} {i + j \choose 2i
- j} x^i y^j.
\end{equation*}
\end{exa}
\begin{proof}

We have
\begin{align*}
\sum_{i,j\ge 0} {i+j\choose 2i-j} x^i y^j &=
 \sum_{i,j\ge 0} \left( \ct (1+t)^{2i-j} (1+\tt)^{2j-i}\right)
x^i y^j \\
&=\ct \sum_{i,j\ge 0} \frac{(1+t)^{2i}}{(1+\tt)^i} x^i \cdot
\frac{(1+\tt)^{2j}}{(1+t)^j} y^j \\
&= \ct \frac{1}{(1-tx-t^2x)(1-\tt y-t^{-2} y)}
\end{align*}
Using Maple, we find the partial fraction expansion in $t$:
\begin{multline}\label{e-exa1} \frac{1}{(1-tx-t^2x)(1-\tt y-t^{-2} y)} =
\frac{1}{1-3xy-x^2y-xy^2}\\
\left(\frac{1-xy+txy+tx^2y}{1-tx-t^2x}
-\frac{y(1+t+tx-xy)}{t^2-ty-y}\right) 
\end{multline}

It is easy to see that
$$\ct \frac{1}{1-3xy-x^2y-xy^2}\frac{1-xy+txy+tx^2y}{1-tx-t^2x} =
\frac{1-xy}{1-3xy-x^2y-xy^2},$$ obtained by setting $t=0$, since
it is a formal power series in $x$ and $y$, with coefficients in
$\mathbb{Q}[t]$. Similarly,
$$\ct \frac{y(1+t+tx-xy)}{t^2-ty-y} =\ct \tt\frac{y(1+\tt+x-\tt
xy)}{1-\tt y-t^{-2} y}=0 ,$$ since it is a formal power series in
$x$ and $y$ with coefficients in $t^{-1}\mathbb{Q}[t^{-1}]$.

Equation \eqref{e-gf1} then follows.
\end{proof}

Similarly, we can compute the generating function of
${i+j+r\choose 2i-j}$ for nonnegative integer $r$, from which it
is easy to deduce \eqref{e-r}.
\begin{exa}\label{ex-2}
\begin{align}\label{e-gfzzz}
\sum_{i,j,r\ge 0} {i+j+r \choose 2i-j} x^i y^j
z^r=\frac{1-z-xy(1-2z-xz)}{(1-3xy-x^2y-xy^2)(1-2z+z^2-xz)}.
\end{align}
\end{exa}
\begin{proof}
We have
\begin{align*}
\sum_{i,j,r\ge 0} {i+j+r \choose 2i-j} x^i y^j z^r &=
\sum_{i,j,r\ge 0} \left( \ct (1+t)^{2i-j} (1+\tt)^{2j-i+u}\right)
x^i y^j z^r\\
&=\ct \sum_{i,j,r\ge 0} \frac{(1+t)^{2i}}{(1+\tt)^i} x^i \cdot
\frac{(1+\tt)^{2j}}{(1+t)^j} y^j\cdot (1+\tt)^rz^r \\
& = \ct \frac{1}{1-tx-tx^2}\cdot \frac{1}{1-\tt y-t^{-2} y} \cdot
\frac{1}{1-z-\tt z}
\end{align*}
By \eqref{e-exa1}, this equals
$$\ct \frac{1}{1-3xy-x^2y-xy^2}\left(\frac{1-xy+txy+tx^2y}{1-tx-t^2x}
-\frac{y(1+t+tx-xy)}{t^2-ty-y}\right)\frac{1}{1-z-\tt z}.$$ Since
${1}/{(1-z-\tt z)}$ is a formal power series in $z$ with
coefficients in $\mathbb{Q}[t^{-1}]$, we can discard the second
part of \eqref{e-exa1} in our computation. So we have
$$\ct \frac{1}{1-3xy-x^2y-xy^2}\left(\frac{1-xy+txy+tx^2y}{1-tx-t^2x}
\right)\frac{1}{1-z-\tt z}.$$ Converting this into partial
fraction in $t$, we get two parts, one with coefficients in
$\mathbb{Q}[t]$, the other with coefficients in
$t^{-1}\mathbb{Q}[t^{-1}]$. So we can discard the second part and
set $t=0$ to get equation \eqref{e-gfzzz}.
\end{proof}

\begin{proof}[Algebraic Proof of Lemma \ref{l-mainlemma}]

For $s$ with $1\le s \le r$ and any $P\in \pa (i,j)$, we factor
$P$, according to its intersection with the line $y=0$, uniquely
as $P_1P_2$, where $P_1$ is a $T^{(r)}$-path, and
 $P_2$ is an $H_r$ path starting with a vertical step. Then
$S(P_2)=(kr+s,0)$, for some $0\le k\le j$. See the left picture of
Figure \ref{fig-9}.

\begin{figure}[ht]
\begin{center}
\input{pic9.pstex_t}
\end{center}
\caption{Example paths for $\pa(i,j)$ and
$\pb(i,j)$.\label{fig-9}}
\end{figure}
For any $k$, there are $\rr{T}(kr+s,ir-s)$ choices for $P_1$ since
it is a \trpath from $(0,-(ir-s))$ to $(kr+s,0)$. There is only
one choice for $P_2$ since it is  a $(0,1)$ step followed by a
fixed number of $(r,0)$ steps. Conversely, any such $P_1P_2$ stays
below the diagonal and hence belongs to $\pa(i,j)$. Summing on all
possible $k$, we get a formula for $|\pa(i,j)|$:
\begin{align}
|\pa(i,j)| &=  \sum_{k=0}^{j} \rr{T}(kr+s,ir-s). \label{e-t-v}
\end{align}

Multiplying both sides of equation \eqref{e-t-v} by $x^iy^j$, and
summing
 on $i$
and $j$, we have
\begin{align*}
\sum_{i,j=0}^\infty |\pa(i,j)|x^iy^j &= \sum_{i,j=0}^\infty
\sum_{k=1}^i \rr{T}(kr+s,ir-s) x^iy^j
\\ &=\sum_{i,j=0}^\infty\sum_{k=0}^i (\ct t^s \alpha^k\beta^i)x^iy^j \\
&=\ct t^s\sum_{j=0}^\infty \beta^ix^i \sum_{k=0}^\infty
\alpha^ky^k\sum_{j\ge
k}y^{j-k}\\
 &=\ct t^s\frac1{(1-\alpha y)(1-\beta x)(1-y)}
\end{align*}

Similarly, any $Q\in \pb(i,j)$ can be factored, according to its
intersection with the line $x=0$, uniquely as $Q_1Q_2$, where
$Q_1$ is a $V_r$ path ending with a horizontal step and $Q_2$ is a
$T^{(r)}$-path. Then $E(Q_1)=(0,-(k'r-s+1))$ for  some $1\le k'
\le i$. See the right picture of Figure \ref{fig-9}. We see that
even in the case $s=1$, $k'$ cannot be zero, since otherwise $Q_1$
will go above the diagonal.

For any $k'$, there is only one choice for $Q_1$ since it is a
fixed number of $(0,r)$ steps followed by a $(1,0)$ step. There
are $\rr{T}(jr+s-1, k'r-s+1)$ choices for $Q_2$ since it is a
\trpath from $(0,-(k'r-s+1))$ to $(jr+s-1,0)$. Summing on all the
possible $k'$, we get a formula for $|\pb(i,j)|$.
\begin{align}
|\pb(i,j)| &=\sum_{k'=1}^i \rr{T}(jr+s-1, k'r-s+1). \label{e-t-h}
\end{align}

A similar computation shows that
\begin{align*}
\sum_{i,j=0}^\infty |\pb(i,j)|x^iy^j &= \ct
t^{s-1}\sum_{j=0}^\infty \alpha^jy^j \sum_{k'=1}^\infty
\beta^{k'}x^{k'}\sum_{i\ge k'}x^{i-k'}\\
&=\ct t^{s-1}\frac{1}{1-\alpha y}\left(\frac{1}{1-\beta
x}-1\right)\frac{1}{1-x}\\
&=\ct t^{s-1}\frac{1}{(1-\alpha y)(1-\beta x)(1-x)}- \ct
t^{s-1}\frac{1}{(1-\alpha y)(1-x)} \\
&=\ct t^{s-1}\frac{1}{(1-\alpha y)(1-\beta x)(1-x)}
-\delta_{s,1}\frac{1}{(1-y)(1-x)}.
\end{align*}
To compute the two generating functions does not seem easy, but
their difference has a simple form.
\begin{multline*}
\ct t^s\frac1{1-\alpha y}\frac{1}{1-\beta x}\frac{1}{1-x}- \ct
t^{s-1}\frac{1}{1-\alpha y}\frac{1}{1-\beta x}\frac{1}{1-y} \\=\ct
t^{s-1} \frac{t-ty-1+x}{(1-x)(1-y)(1-\alpha y)(1-\beta x)}
\end{multline*}
Direct computation shows that
\begin{align*}
\frac{t-ty-1+x}{(1-\alpha y)(1-\beta x)} &=- \frac{x}{t^r(1-\beta
x)}-\frac{1-y-t}{1-\alpha y}
\end{align*}
It is easy to see that $-\displaystyle\frac{1-y-t}{1-\alpha y}$,
belonging to $\QQ[t][[y]]$, has constant term constant term $-1$
in $t$, and $-\displaystyle\frac{x}{t^r(1-\beta x)}$, belonging to
$t^{-1}\QQ[t^{-1}][[x]]$, has constant term $0$. Hence
$$\ct t^{s-1} \frac{t-ty-1+x}{(1-x)(1-y)(1-\alpha y)(1-\beta x)} =-\delta_{s,1}
\frac{1}{(1-x)(1-y)}.$$

Put the above altogether, we obtain
\[\sum_{i,j=0}^\infty (| \pa(i,j)|-|\pb(i,j)|) x^iy^j =0.
\qedhere\]
\end{proof}

\begin{lem}\label{l-lemma}
\begin{align}\label{e-gf-tr0}
\sum_{i,j\ge 0} \rr{T}(ri,rj) x^i y^j =
 \frac{ x(1-x)^{r-1}-y(1-y)^{r-1}}{x(1-x)^r-y(1-y)^r}.
\end{align}
\end{lem}
\begin{proof}
Let $x= uv$ and $y= v$. Then
\begin{align*}
\sum_{i,j\ge 0} \rr{T}(ri,rj) x^i y^j &= \sum_{i,j\ge 0}
\rr{T}(ri,rj) u^i
v^{i+j}\\ &=\sum_{n\ge i\ge 0} \rr{T}(ri,r(n-i)) u^i v^n \\
&=\sum_{n\ge i\ge 0} \ct \frac{\alpha^n}{t^{ri}} u^i v^n \\
&= \ct \frac{1}{1-ut^{-r}} \frac{1}{1-\alpha v}.
\end{align*}
We have
\begin{align*}
\frac{1}{1-\alpha v} &= \frac{1-t}{(1-t)-v(1-t^{r+1})} \\
&= \frac1{1-v} \frac{1-t}{1-t \displaystyle\frac{1-t^rv}{1-v}}\\
&= \frac1{1-v}
\frac{(1-t)\left(1+t\displaystyle\frac{1-t^rv}{1-v}+\cdots +
\left(t\displaystyle\frac{1-t^rv}{1-v}\right)^{r-1}\right)}
{1-t^r\left(t\displaystyle\frac{1-t^rv}{1-v}
\right)^r}\\
&= \frac{1}{1-v}
\frac{1-t^r\left(\displaystyle\frac{1-t^rv}{1-v}\right)^{r-1}
}{1-t^r\left(t\displaystyle\frac{1-t^rv}{1-v}\right)^r}
+\mbox{other terms}
\end{align*}
Since the other terms contain only terms like $a_{rm+s} t^{rm+s}$
for  $1\le s\le r-1$, they do not contribute to the constant term
in $t$. Let $z=t^r$. Then we have
\begin{align*}
\ct \frac{1}{1-u t^{-r}}\frac{1}{1-\alpha v} &= \ct _z\frac{1}{1-u
z^{-1}}\frac{1}{1-v}
\frac{1-z\left(\displaystyle\frac{1-zv}{1-v}\right)^{r-1}}{1-z\left(\displaystyle\frac{1-zv}{1-v}\right)^r
},
\end{align*}
where $\ct _z$ means to take the constant term of a function in
$z$. Since the other part of the right side of the above equation
is a formal power series in $z$, it is straightforward to obtain
$$
\ct \frac{1}{1-u t^{-r}}\frac{1}{1-\alpha v} = \frac{1}{1-v}
\frac{1-u
\left(\displaystyle\frac{1-uv}{1-v}\right)^{r-1}}{1-u\left(\displaystyle\frac{1-uv}{1-v}\right)^r}.
$$
Replacing $u$ with $x/y$ and $v$ with $y$, we get formula
\eqref{e-gf-tr0}
\end{proof}

\begin{thm}
$$\det\left(\rr{g}_{i+j}\right)\ijn =\det \left(\rr{T}(ri,rj)\right)\ijn =\det \left(T^{(r-1)}(rj-i,ri-j)\right)\ijn
.$$
\end{thm}
Note that the identities in this theorem appeared in \eqref{e-98}
and \eqref{e-99}.
\begin{proof}
We use the technique of Section 4. The generating function for the
first determinant is
$$\frac{x \rr{g}(x)-y\rr{g}(y)}{x-y}.$$
Since $\rr{f}(x)=\rr{g}(x)-1$ is a formal power series in $x$
satisfying
$$\rr{f}(x)=x/\left(1+\rr{f}(x)\right)^{r+1},$$
we make
the substitution $x\to x/(1+x)^{r+1}$ and $y\to y/(1+y)^{r+1}$.
The generating function becomes
$$\frac{x(1+x)^{-r-1}(1+x)-y(1+y)^{-r-1}(1+y)}{x(1+x)^{-r-1}- y(1+y)^{-r-1}}
=\frac{x(1+x)^{-r}-y(1+y)^{-r}}{x(1+x)^{-r-1}- y(1+y)^{-r-1}}.$$
Normalizing and dividing by $(1+x)(1+y)$, we get
\begin{align}
\label{e-th9.5} \frac{
x(1-x)^{r-1}-y(1-y)^{r-1}}{x(1-x)^r-y(1-y)^r}.
\end{align}
The first equality hence follows from Lemma \ref{l-lemma}.

To show the second equality, we compute the generating function
for the third determinant.
\begin{align}
&\sum_{i,j\ge 0} T^{(r-1)}(rj-i,ri-j)x^iy^j \nonumber\\
&\qquad \qquad= \sum_{i,j\ge 0}
T^{(r-1)}((r-1)j+(j-i),(r-1)i-(j-i))x^iy^j \nonumber \\
&\qquad \qquad=\sum_{i,j\ge 0} \ct t^{j-i}(1+t+\cdots +t^{r-1})^j
(1+\tt+\cdots
+t^{-r+1})^i x^i y^j \nonumber \\
&\qquad \qquad=\ct \sum_{i,j\ge 0} (\tt+t^{-2}\cdots
+t^{-r})^i x^i (t+t^2+\cdots +t^{r})^jy^j \nonumber \\
&\qquad \qquad=\ct \frac{1}{1-(\tt+t^{-2}\cdots +t^{-r}) x}
\frac{1}{1-(t+t^2+\cdots +t^{r})y}\label{e-s9-third}
\end{align}
Similarly, the generating function for the second determinant is
\begin{align} \ct \frac{1}{1-(1+\tt+t^{-2}\cdots +t^{-r}) x}
\frac{1}{1-(1+t+t^2+\cdots +t^{r})y}.\label{e-s9-second}
\end{align}

The following computation shows that the \eqref{e-s9-second} can
be obtained from \eqref{e-s9-third} by making the substitution
$x\to x/(1-x)$ and $y\to y/(1-y)$, and then dividing by
$(1-x)(1-y)$. This yields the second equality.
\begin{align*}
&\ct \frac{1}{1-(1+\tt+t^{-2}\cdots
+t^{-r}) x} \frac{1}{1-(1+t+t^2+\cdots +t^{r})y} \\
&= \ct \frac{1}{1-x-(\tt+t^{-2}\cdots
+t^{-r}) x} \frac{1}{1-y-(t+t^2+\cdots +t^{r})y}\\
&=\ct \frac{1}{(1-x)(1-y)} \frac{1}{1-(\tt+t^{-2}\cdots +t^{-r})
\frac{x}{1-x}} \frac{1}{1-(t+t^2+\cdots
+t^{r})\frac{y}{1-y}}.\qedhere
\end{align*}

\end{proof}

{\bf Acknowledgements:} The authors wish to thank the
Institut Mittag-Leffler and the organizers of the
Algebraic Combinatorics program held there in the spring of 2005,
Anders Bj\"{o}rner and Richard Stanley.
\bibliographystyle{amsplain}

\providecommand{\bysame}{\leavevmode\hbox
to3em{\hrulefill}\thinspace}
\providecommand{\MR}{\relax\ifhmode\unskip\space\fi MR }
\providecommand{\MRhref}[2]{%
  \href{http://www.ams.org/mathscinet-getitem?mr=#1}{#2}
} \providecommand{\href}[2]{#2}

\end{document}